\documentclass[11pt]{article}

\usepackage[margin=1in]{geometry}
\usepackage{cite}
\usepackage{graphics} 
\usepackage{graphicx}
\usepackage{epsfig}
\usepackage{tikz,pmat}
\usepackage{enumitem}
\usepackage{gensymb}

\usepackage{amsmath} 
\usepackage{amssymb}  
\usepackage{url}
\usepackage{algpseudocode,algorithm,algorithmicx}
\usepackage[subnum]{cases}
\usepackage{empheq}
\PassOptionsToPackage{hyphens}{url}
\usepackage{array}
\newcolumntype{L}[1]{>{\raggedright\let\newline\\\arraybackslash\hspace{0pt}}m{#1}}
\newcolumntype{C}[1]{>{\centering\let\newline\\\arraybackslash\hspace{0pt}}m{#1}}
\newcolumntype{R}[1]{>{\raggedleft\let\newline\\\arraybackslash\hspace{0pt}}m{#1}}

\graphicspath{{./Figures/}}

\newcommand{\x}{\mathbf{x}}

\newcommand{\uu}{\mathbf{u}}
\newcommand{\yy}{\mathbf{y}}
\newcommand{\f}{\mathbf{f}}

\newcommand{\A}{\mathbf{A}}
\newcommand{\B}{\mathbf{B}}
\newcommand{\C}{\mathbf{C}}

\newcommand{\X}{\mathcal{X}}
\newcommand{\G}{\,\mathcal{G}}
\newcommand{\U}{\mathcal{U}}
\newcommand{\E}{\mathcal{E}}
\newcommand{\R}{\mathbb{R}}
\newcommand{\ETA}{\pmb{\eta}}

\newtheorem{theorem}{Theorem}
\newtheorem{definition}{Definition}
\newtheorem{assumption}{Assumption}
\newtheorem{lemma}{Lemma}
\newtheorem{remark}{Remark}
\newtheorem{proposition}{Proposition}
\newtheorem{coro}{Corollary}

\newtheorem{problem}{Problem}
\newenvironment{proof}{\paragraph*{Proof:}}{\hfill$\square$}

\mathchardef\mhyphen="2D

\def\beq{\begin{equation}}

\def\eeq{\end{equation}}

\title{Guaranteed Model-Based Fault Detection in Cyber-Physical Systems: A Model Invalidation Approach}

\author{Farshad Harirchi and Necmiye Ozay
\thanks{The authors are with the Electrical Engineering and Computer Science Department, University of Michigan, Ann Arbor, MI 48109. {\tt{\{harirchi,necmiye\}@umich.edu}}}}%

\begin{document}
\maketitle

\begin{abstract}
This paper presents a sound and complete fault detection approach for cyber-physical systems represented 
by hidden-mode switched affine models with time varying parametric uncertainty. The fault detection approach builds upon 
techniques from model invalidation. In particular, a set-membership approach is taken where the noisy 
input-output data is compared to the set of behaviors of a nominal model. As we show, this set-membership 
check can be reduced to the feasibility of a mixed-integer linear programming (MILP) problem, which can be 
solved efficiently by leveraging the state-of-the-art MILP solvers. In the second part of the paper, given a system 
model and a fault model, the concept of $T$-detectability is introduced. If a pair of system and fault models satisfies 
$T$-detectability property for a finite $T$, this allows the model invalidation algorithm to be implemented in a 
receding horizon manner, without compromising detection guarantees. 
In addition, the concept of weak-detectability is introduced which extends the proposed 
approach to a more expressive class of fault models that capture language constraints on the mode sequences. 
Finally, the efficiency of the approach is illustrated with numerical examples motivated by smart building radiant systems.
\end{abstract}


\section{Introduction}

Cyber-physical systems are combinations of physical processes and embedded computers. 
The embedded computers collect data from the process through sensors and control it in a
closed-loop manner. With the increase in data acquisition and storage capacity and the decrease
in sensor costs, it is possible to collect large amounts of data during the operation of complex
cyber-physical systems. For instance, ``a four-engine jumbo jet can create 640 terabytes of data
in just one crossing of the Atlantic Ocean" \cite{Rajah2014Taming}. As discussed in 
\cite{Sznaier2014Surviving}, this exponential growth in the data collection capabilities 
is a major challenge for systems and control community. Sensor/information-rich networked 
cyber-physical systems, from air traffic or energy networks to smart buildings, are getting tightly
integrated into our daily lives. As such, their safety-criticality increases. For such systems, it is 
crucial to detect faults or anomalies in real-time to support the decision-making process and 
to prevent potential large-scale failures.

\subsection{Contributions}

This paper presents a fault detection scheme to enhance the reliability of a class of cyber-physical
systems that are represented by hidden-mode switched affine models with time-varying parametric
uncertainty. This modeling framework is quite expressive and can be used to describe a wide range
of cyber-physical systems such as heating, ventilation and air conditioning (HVAC) systems in smart
buildings \cite{Weimer2013Parameter}, wind turbines \cite{Burkart2011Nonlinear}, 
power systems and power electronics, automotive systems, aircraft, air traffic, and network and 
congestion \cite{Sun2006Switched}. We model faults also in this framework allowing us to capture 
many scenarios including cascaded faults or various types of cyber or physical attacks. Note that 
linear time invariant systems with or without noise or affine parametric uncertainty are special cases
of the modeling framework. 

The proposed fault detection scheme builds on set-membership model
invalidation approaches \cite{Rosa2010Fault,Ozay2010Model,OzayConvex2014}.
Unlike many set-membership methods that compute explicit
set representations and propagate them via set-valued observers, the proposed
method uses an online optimization formulation 
for model invalidation that keeps implicit constraints to represent sets.
In particular, we show that model invalidation problem for this
class of systems can be reduced to the feasibility of a 
MILP problem, which can be 
checked efficiently using state-of-the-art solvers \cite{cplex}. 
Additionally, the concept of $T$-detectability is introduced, which enables us to apply the 
model invalidation approach for fault detection in a receding horizon manner (with a horizon size $T$) 
without losing detection guarantees. 
Even though there are some practical systems with fault models that are $T$-detectable, 
a limited class of system and fault models satisfy this property. We further
discuss weak detectability that incorporates language constraints on the switching
sequences of the faults to enable detection of a broader class of faults.
Algorithms that can be used to find the minimum $T$ (if it exists) are presented.

A preliminary version of this paper is published in \cite{Harirchi2015Model}. 
The current paper significantly extends the model class by
allowing uncertainty in variables and considering language constraints in
the mode sequence. Moreover, more efficient MILP-based necessary and sufficient
conditions for verifying $T$-detectability are given and some
connections to mode observability of switched systems are pointed out.

\subsection{Literature Review on Fault Detection}
Model-based fault detection has a long history
starting with early work on failure detection filters \cite{Beard1971Failure,Jones1973Failure}. 
A vast majority of fault detection methods are based on residual generation,
where the residual is evaluated by simple thresholding methods 
or more complicated classifiers to decide between faulty and normal behaviors 
\cite{Simani2003Model,Isermann2006Fault,Ding2008Model,Patton2013Issues}. 
The residual generation
 methods are classified into three main categories \cite{Simani2003Model}: 
 parameter estimation-based \cite{Isermann1993Fault}, 
 observer-based \cite{Frank1993Advances,Patton1997Observer,Shames2011Distributed}
  and parity equation-based \cite{Gertler1997Fault} techniques. 
  All these residual generation techniques can be implemented in real-time, 
  but even when a specific fault model exists, 
their behavior is usually analyzed only asymptotically and they fail to 
provide any finite-time detection guarantees. 

As an alternative to residual generation,
set-membership fault detection methods have been proposed
both for passive \cite{Rosa2010Fault, Rosa2013Fault} 
and active \cite{Nikoukhah1998Guaranteed,Nikoukhah2006Auxiliary,Scott2014Input} fault detection.
Some of these methods proceed by computing convex-hulls of potentially non-convex reachable sets of the system 
and comparing the actual output to this set \cite{Rosa2010Fault}. Since, reachable sets are over-approximated, 
this leads to only sufficient conditions, that is, they guarantee there are no false alarms. 
Scott {\em et al.} \cite{Scott2014Input} pose a mixed-integer quadratic programming problem
to find an optimal separating input to detect faults. They use zonotopes to represent 
and propagate state constraints. The order of the zonotopes is used to trade-off between
the complexity and the conservativeness of the approach. This approach is extended in \cite{Raimondo2016Closed} to handle admissible sets in the 
form of constrained zonotopes and to take advantage of real-time measurements to calculate an improved separating input. 
Despite the fact that all the above mentioned set-membership methods are proposed for linear systems, their scalability
is somewhat limited to be applied in real-time \cite{Scott2014Input}.

Non-linear and hybrid systems have also attracted notable attention from fault detection community. 
Most of the research is concentrated around residual generation type methods \cite{Garcia1997Deterministic, Hammouri1999Observer, Pan2015Online}. 
De Persis and Isidori \cite{De2001Geometric} develop analytical necessary and sufficient conditions under which the problem
of fault detection and isolation becomes solvable for non-linear systems. However, they do not particularly address the computational 
aspects of the problem. 
Observer-based methods are employed for fault diagnosis in hybrid systems
both when the discrete mode is observed  \cite{Mcilraith2000Hybrid} or hidden \cite{Narasimhan2007Model},
using variants of Kalman filters.
In recent work \cite{Deng2015Verification}, Deng {\em et al.} investigate fault diagnosis problem in hybrid systems, by constructing a finite 
abstraction for the hybrid automaton and analyzing the diagnosability of the abstract system using tools from 
discrete-event systems, providing some detection guarantees.
However, the mode signal is assumed to be observed in \cite{Deng2015Verification}. 

The notion of $T$-detectability is closely related to distinguishability of dynamical systems \cite{Grewal1976Identifiability,De2016Observability}
and observability in switched systems \cite{vidal2002observability,Babaali2004Observability}.
If two systems are distinguishable, then there exists a non-zero input which makes their trajectories different at least for one time instance. In \cite{De2011Location}, it is illustrated that if two models are distinguishable for a particular control input, then they are distinguishable for a generic control law. Lou {\em et. al} \cite{Lou2009Distinguishability} introduced the concept of input-distinguishability, which restricts distinguishability to all non-zero inputs contained in a convex set. Rosa and Silvestre \cite{Rosa2011Distinguishability} consider the input-distinguishability for discrete time linear time-invariant systems with bounded disturbance and noise. They provide a necessary and sufficient rank condition to check for a given size of time horizon $(T)$ if two systems are input distinguishable or not. This condition is seldom satisfied in practice, hence they added extra constraint which enforces the persistence of excitation for disturbance \cite{Rosa2011Distinguishability}. The concept of $T$-detectability introduced here is closely related to absolute input-distinguishability \cite{Rosa2011Distinguishability}. In other words, $T$-detectability of a pair of hidden-mode switched affine models subject to process and measurement noise provides necessary and sufficient conditions for their input-distinguishability. The time horizon $T$ in detectability is an upper bound on the detection delays (time from the occurrence of fault to detection alarm) \cite{Adnan2011Expected,Mariton1989Detection,Stoorvogel2001Delays}. 

Model invalidation was originally proposed as a way to build trust 
in models obtained through a system identification step or discard/improve them before using these models 
in robust control design \cite{smith1992model}.  In model invalidation problem, one starts with
a family of models (i.e., a priori or admissible model set) and experimental input-output data collected from
a system (i.e., a finite execution trace) and tries to determine whether the experimental data can be
generated by one of the models in the initial model family. Its relation to fault and anomaly detection
  has been pointed out in \cite{Rosa2010Fault,Ozay2010Model,Cheng2012Convex,OzayConvex2014}
   for linear time-varying systems and hybrid systems in autoregressive form. 
   A fault detection scheme based on model-invalidation 
   for polynomial state space models subject to noise and uncertainty in parameters is 
   recently proposed in \cite{Harirchi2016Model}. The convex relaxations
   for model invalidation problem mostly provide sufficient conditions that
   can be efficiently checked to detect faults but only necessary for large relaxation orders.

\noindent {\emph{Notation:}} Let $ \x \in \R^n$ denote a vector and $\x(i)$ indicate its $i^{th}$ element. Also, let $ \A \in \R^{n\times m}$ represent 
a matrix and $\A(i,j)$ indicate the element on the $i^{th}$ row and $j^{th}$ column of the matrix $\A$. The row range and null space of a matrix $\A$ 
are denoted by $\text{Range}(\A)$ and $\text{Null}(\A)$, respectively. The infinity norm 
of a vector $\textbf{x}$ is given by $\| \textbf{x} \|  \doteq \max_i \x(i)$. The set of positive and non-negative integers 
up to $n$ are denoted by $\mathbb{Z}_n^{+} $ and $\mathbb{Z}_n^{0}$, respectively. A vector of all ones of appropriate size is denoted by $\mathbf{1}$. The 
Hadamard product of two matrices $\A_1$ and $\A_2$ of the same dimensions is indicated by $\A = \A_1\odot\A_2$, 
and defined as $\A(i,j) \doteq \A_1(i,j)\A_2(i,j)$ for all $i,j$. 

\section{Modeling Framework} \label{sec:model}
In this section, we describe the class of systems and the modeling framework considered in this paper.  
In particular, we consider discrete-time hidden-mode switched affine models 
with time-varying parametric uncertainty that are subject to noise. 
We refer to these models as SWA models. SWA models consist of a collection of affine models,
which we define first.

\begin{definition}\label{def:UA}
(Affine Model) An affine model $G^\Delta$ with time-varying parametric uncertainty 
has the following form:
\begin{equation}\label{eq:affine_mod}
\begin{aligned}
\x_{k+1}  	& = (\A+\A^\Delta_k) \x_k +  (\B+ \B^\Delta_k)   \uu_k + \f +\f^\Delta_k\\
\yy_k  	& = (\C+\C^\Delta_{k})  \x_k +\ETA_k
\end{aligned},
\end{equation}
where $\x_k, \uu_k, \yy_k$ are the state, input and output at time $k$,
$\A,\B, \C, \f$ are the system matrices; 
and $\A^\Delta_k,\B^\Delta_k, \C^\Delta_k, \f^\Delta_k$
are the weighted uncertain variable matrices affecting the system at time $k$.
In particular, for each $\pmb{\tau}^{\Delta}_k, \; \pmb{\tau} \in \{ \A,\B, \C, \f\}$,
we assume:
\begin{equation} \label{eqn:uncnormal}
\pmb{\tau}^{\Delta}_k = \hat{\pmb{\tau}} \odot \Delta^{\tau}_{k}
\end{equation} 
for some normalization matrix $\hat{{\pmb{\tau}}}$ with nonnegative entries and uncertainty matrix $\Delta^{\tau}_{k}$
whose entries satisfy
\begin{equation} \label{eqn:unccompwise}
-1 \leq \Delta^{\tau}_{k}(m,l)\leq 1,
\end{equation}
where $\Delta^{\tau}_{k}(m,l)$ corresponds to the uncertainty variable associated with the term 
on $m^{th}$ row and $l^{th}$ column of parameter $\tau$ at time $k$.
\end{definition}

For simplicity of notation, we collect the entries of $\Delta^{\tau}_{k}$ for all $\pmb{\tau} \in \{ \A,\B, \C, \f\}$, 
in a vector $\Delta_k \in \Omega \subset \R^{n_\Delta}$, where 
$n_{\Delta}$ is the total number of uncertain parameters. 
It follows from the assumption \eqref{eqn:unccompwise} that
the set $\Omega$ of admissible uncertainties is the unit infinity-norm ball in $\R^{n_\Delta}$.

Now, we define SWA models, the main modeling framework used in this paper.
\begin{definition}
	(Switched Affine Model) A switched affine (SWA) model is defined by:
	\begin{equation} \label{eqn:USWA}
	\mathcal{G} = (\X, \E, \U, \{ G_i^\Delta\}_{i=1}^s),
	\end{equation} 
	where $\X \subset \R^n$ is the set of states, $\E \subset \R^{n_y}$ is the set of measurement noise values, 
	$\mathcal{U} \subset \R^{n_u}$ is the set of inputs. The collection $\{G^{\Delta}_i\}_{i=1}^s$ 
	is the set of $s$ affine models. The evolution of $\G$ is governed by:
\begin{equation} \label{eqn:SIMUSWA}
\begin{aligned} 
\x_{k+1}  	& = (\A_{\sigma_k}+\A^\Delta_{\sigma_k,k}) \x_k +  (\B_{\sigma_k}+\B^\Delta_{\sigma_k,k})   \mathbf{u}_k+ \\ 
& \quad \quad \mathbf{f}_{\sigma_k} +\f^\Delta_{\sigma_k,k} \\
\mathbf{y}_k  	& = (\C_{\sigma_k}+\C^\Delta_{\sigma_k,k})  \x_k +\ETA_k 
\end{aligned},
\end{equation}
where $\sigma_k \in \{1, \hdots , s \}$ indicates the active affine model, i.e., the mode, at time $k$, 
and $\A_i,\B_i, \C_i, \f_i$ are the system matrices of the affine model $G_i^\Delta$.
\end{definition}
In what follows, we assume that the sets of admissible states, inputs and noise are hyper-rectangles
defined as:  
\begin{align}
		\mathcal{X} &=\{ \x \mid X_l \leq \x \leq X_u \} \label{eqn:stateconst}, \\
		\mathcal{E} &= \{\ETA \mid  \epsilon_l \leq \ETA \leq \epsilon_u \} \label{eqn:noiseconst},\\
		\U &= \{\uu \mid U_{l} \leq \uu \leq U_u \}, \label{eqn:inputbounds}
\end{align}
where subscript $l$ and $u$ denote lower and upper bounds, respectively.
The bound on the state usually represents the part of the state-space where
the model is valid and can be taken to be infinite. The distinction between the noise $\ETA$ and 
the inputs $\uu$ is that the former is not measured at run-time but the latter is measured
and can capture operating conditions or effects of other subsystems on the specific system under
consideration. Therefore, typically both the noise and inputs are bounded in some range.
If a closed-loop system is considered and the feedback law is already in place,
it is possible to omit the inputs, and the proposed methods are trivially applicable in that case as well.
	
For the sake of simplicity in notation, we omit process noise, possible affine term in the output equation 
and the feed-forward term in the output equation in Def.~\ref{def:UA}. 
It is, however, straightforward to apply the proposed techniques 
when considering all the mentioned terms. 
Similarly, the uncertainty model can also be generalized to capture some correlations between entries. 


\section{Problem Statements} \label{sec:problem}
In this section, we present some preliminary definitions together with the two problems,
solutions of which provide the basis of the proposed fault detection approach.

\begin{definition}
	(length-$N$ behavior) The \emph{length-$N$ behavior} associated with an SWA model $\G$ is the set of all 
	length-$N$ input-output trajectories compatible with $\G$, given by the set
	\begin{align*}
	\mathcal{B}^N_{swa}(\G) \doteq \big \{ \{\uu_k, \yy_k\}_{k=0}^{N-1} \mid \uu_k \in \mathcal{U},  \; \exists \x_
	k\in \mathcal{X},\ETA_k \in \mathcal{E},\\ 
	\Delta_k \in \Omega, \sigma_k\in \mathbb{Z}_s^+: \eqref{eqn:SIMUSWA} 
	\text{ holds for all } k\in\mathbb{Z}_{N-1}^0  \big \}.
	\end{align*}
	With slight abuse of terminology, 
	when $N$ is clear from the context, 
	we call $\mathcal{B}^N_{swa}(\G)$ just the \emph{behavior} of $\G$.
\end{definition} 

The consistency set associated with a behavior of an SWA model is defined as follows.
\begin{definition}
	(consistency set) Let $\{\uu_k,\yy_k \}_{k=0}^{N}$ be an input-output trajectory over a time window $[0,N]$. 
	The {\em consistency set} associated with $\{\uu_k,\yy_k \}_{k=0}^{N}$ of an SWA model $\G$ is defined as follows:
	\begin{align} \label{eqn:consistency}
	&\mathcal{T}_{\G}\big ( \{\uu_k, \yy_k \}_{k=0}^{N}\big ) \doteq   \big \{ \x_{0:N},\ETA_{0:N}, \Delta_{0:N}, \sigma_{0:N}
	 \mid  \nonumber \\ 
	& \quad \quad \quad \eqref{eqn:SIMUSWA}\mhyphen\eqref{eqn:noiseconst} \text{ hold for all } k \in \mathbb{Z}_N^0 \big \}, 
	\end{align}	
	where subscript $0:N$ indicates all the samples between time instances $0$ and $N$.
\end{definition}
In words, the consistency set is the set of state, noise, uncertainty and mode sequences that evolves
 through dynamics of the SWA model \eqref{eqn:SIMUSWA}, and given the input sequence $\{\uu_k \}_{k=0}^{N} \in \mathcal{U}^{N+1}$ 
 generates the output sequence $\{\yy_k \}_{k=0}^{N}$. As the input sequence is needed to calculate the consistency set, 
 we can check whether the bounds on $u_k$, (cf. \eqref{eqn:inputbounds}) are violated or not, separately; 
 and if they are violated we take $\mathcal{T}_{\G}\big ( \{\uu_k, \yy_k \}_{k=0}^{N}\big )=\emptyset$ trivially.
 
 The first problem we address in this paper is the model invalidation problem. 
 Roughly speaking, given an input-output trajectory and an SWA model, 
 the model invalidation problem is to determine whether or not the data is compatible 
with the model. This problem can be formally stated in terms of behaviors or 
consistency sets as follows.

\begin{problem} \label{Pro:MIP}
	Given $ \{\uu_k, \yy_k \}_{k=0}^{N-1}$, an input-output trajectory, and an SWA model $\G$,
	 determine whether or not the input-output trajectory is contained in the behavior of $\G$. That 
	 is, determine whether or not the following is true
	\begin{align} \label{eqn:MIPDEF}
	 \{\uu_k,\yy_k \}_{k=0}^{N-1} \in \mathcal{B}_{swa}^N(\G),
	\end{align}
	or equivalently 
	$\mathcal{T}_{\G} \big (\{\uu_k,\yy_k \}_{k=0}^{N-1} \big ) \neq \emptyset$.
\end{problem}
 
For the well-posedness of the problems studied, 
we assume that $\mathcal{B}^N_{swa}(\G) \neq \emptyset$ 
for any positive integer $N$. This guarantees that the model is capable of
generating some infinite trajectories; therefore, it does not invalidate
itself if one waits long enough. 

 The abnormal trajectories for a system are those that cannot be generated by
the model of a system. 

\begin{definition} (abnormal trajectory) An input-output trajectory $\{\uu_k,{\yy}_k\}_{k=0}^{N-1}$ is called
	 \emph{abnormal} for an SWA model $\G$ if
	\[
	\{\uu_k,{\yy}_k\}_{k=0}^{N-1} \notin \mathcal{B}^N_{swa}(\G),
	\]
	or equivalently $\mathcal{T}_{\G}(\big \{\uu_k, \yy_k \big \}_{k=0}^{N-1}) = \emptyset$.	
\end{definition}

With this definition, it is clear that a trajectory being abnormal is equivalent to the model being invalid
for that trajectory.
Therefore, a solution to Problem~\ref{Pro:MIP} can be readily used to detect abnormal trajectories
or anomalies. 
Note that detecting abnormal trajectories does not require explicit models for the anomaly. 
Given that a cyber-physical system can fail (or be attacked)
 in infinitely many different ways, not needing to model these failure modes is advantageous. 
 On the other hand, 
 if one has an explicit model for a given fault, then this 
 information can be used to develop more efficient fault detection schemes. 
 In this paper, we represent faults also using SWA models. Such fault models can
 be obtained either using first principles in cases when the fault leads to a change in physical
 parameters or using robust system identification techniques in cases where there is
 enough data from the faulty behavior.

\begin{definition}\label{def:fault_mod} (fault) A \emph{fault model} for a system with an SWA model
	$\G= (\X, \E, \U, \{ G_i^\Delta\}_{i=1}^s)$ is another SWA model 
	$\G^f= (\bar{\X}, \bar{\E}, \bar{\U}, \{ \bar{G}_i^{\Delta}\}_{i=1}^{\bar{s}})$ with the same number of states, inputs and outputs. 
\end{definition}

In general, we expect faults to lead to abnormal trajectories. One can argue that a fault that does not lead to
an abnormal trajectory, might not be that important from an operational standpoint. 
Moreover, it is not possible to detect a fault that has trajectories that are identical to the system trajectories.
Next, we define $T$-detectability, which is a property of a given pair of SWA fault and system models, 
that measures how long it takes for a fault to lead to an abnormal trajectory.

\begin{definition}\label{def:T_detec} ($T$-detectability) A fault model $\;\G^f$ for a system model $\G$ is called 
	\emph{$T$-detectable}\footnote{This notion is symmetric and therefore we also say $\G$ and $\G^f$ are $T$-detectable.}, if $\mathcal{B}^T_{swa}(\G) \cap \mathcal{B}^T_{swa}(\G^f) = 
	\emptyset$, where $T$ is a positive integer. 
\end{definition}   

Clearly, if a fault is $T$-detectable for a system, then it is $\bar{T}$-detectable for the same system for all $\bar{T}\geq T$. 

The second problem we are interested in is the
characterization of the $T$-detectability property for a pair of fault and system models.

\begin{problem}
	Given two SWA models, $\G$ and $\G^f$, and an integer $T$, determine
	 whether the fault model $\G^f$ is $T$-detectable for the system model $\G$, or not.
	That is, check if the set
	\begin{equation}
	\mathcal{B}^{T}_{swa}(\G) \cap \mathcal{B}^{T}_{swa}(\G^f), 
	\end{equation}
	is empty or not.
\end{problem}

As we show later in the paper, if a fault model is $T$-detectable for a system model, then 
a solution to the model invalidation problem in a receding horizon manner with a fixed length $T$ can be used to efficiently 
detect faults without compromising the detection guarantees. 


\section{Fault Detection Scheme} \label{sec:approach}
In this section, we propose a fault detection scheme based on model invalidation.
We initially present optimization based solutions to model invalidation and $T$-detectability problems.
 Additionally, we introduce $T$-weak detectability property that allows us to incorporate
constraints on the mode sequences. 
Finally, we present an efficient fault detection scheme for $T$-(weak) detectable faults.

\subsection{Model Invalidation} \label{sec:MI}
As discussed in Section \ref{sec:problem}, given a model of a system,
detecting an abnormal trajectory is equivalent to model being invalid. 
Therefore,
the solution to model invalidation problem, can be readily applied for detecting anomalies in a system
from the input/output measurements.
Next, we define a series of feasibility problems that are equivalent to Problem \ref{Pro:MIP}. 
That is, for a given input-output trajectory $\{\uu_k, \yy_k\}_{k=0}^{N-1}$, 
we encode the consistency set $\mathcal{T}_{\G}\big ( \{\uu_k, \yy_k \}_{k=0}^{N-1}\big )$
with a feasibility problem.
Since we have a hidden-mode switched system model, 
at each time, we require the data to satisfy the dynamic
constraints for at least one mode. We ensure the satisfaction of state and output equations by leveraging ideas from big M formulation of mixed integer programming \cite{Raman1994Modelling, Bemporad1999Control}.
In order to capture this, consider the following mixed integer
non-linear problem:
\begin{subequations} 
	\begin{align} \label{eqn:feasoptnonlineartot} \tag{12}
		\text{Find } & \x_k, \ETA_k, \Delta_k, \; a_{i,k}, \; \forall k \in \mathbb{Z}_{N-1}^0, \forall 
		i\in \mathbb{Z}_s^+ \qquad \qquad \qquad \qquad \qquad \qquad \qquad \qquad \quad \; \; \;  \\
		\text{s.t. } &  \; \forall k \in \mathbb{Z}_{N-1}^0, \forall i\in \mathbb{Z}_s^+: \nonumber
	\end{align}
	\vspace{-0.2in}
	\begin{empheq}[left={ \empheqlbrack}]{align}
		& -\bar{a}_{i,k}M\mathbf{1} \leq \x_{k+1} - (\mathbf{A}_{i}+\A^\Delta_{i,k}) \x_k - (\B_{i}+\B^\Delta_{i,k}) \mathbf{u}_k- \mathbf{f}_{i}- \f^\Delta_{i,k} \leq \bar{a}_{i,k}M\mathbf{1}, \label{eqn:stateconst0} \\
		& -\bar{a}_{i,k}M{\mathbf{1}} \leq \mathbf{y}_k - (\C_{i}+\C^\Delta_{i,k})  \x_k - \ETA_k \leq \bar{a}_{i,k}M{\mathbf{1}},  \label{eqn:stateconst1} \\
		& \textstyle \sum_{i\in \mathbb{Z}_s^+}a_{i,k} = 1, \; a_{i,k} \in \{0,1\}, \label{eqn:stateconstbin} \\
		& X_l \leq \x_k \leq X_u, \; \epsilon_l \leq \ETA_k \leq \epsilon_u, \label{eqn:stateconst2}\\ 
		& \pmb{\tau}^\Delta_{i,k} = \hat{\pmb{\tau}}_i \odot \Delta^{\tau}_{i,k}, \; \pmb{\tau}\in  \{ \A, \B, \C, \f \} \label{eqn:uncertaintyconst}\\
		& |\Delta^{\B}_{i,k}(m,p)| \leq 1, |\Delta^{\f}_{i,k}(m)| \leq 1 \; \forall m\in \mathbb{Z}_n^+, p \in \mathbb{Z}_{n_u}^+\label{eqn:uncertaintyconst2} \\
		& |\Delta^{\A}_{i,k}(m,l)| \leq 1, \; \forall m,l \in \mathbb{Z}_n^+,\label{eqn:uncertaintyconst3} \\
		& |\Delta^{\C}_{i,k}(o,l)| \leq 1, \; \forall o \in \mathbb{Z}_{n_y}^+, \;  \forall l \in \mathbb{Z}_n^+,\label{eqn:uncertaintyconst4}
	\end{empheq}
\end{subequations} 
where $\bar{a}_{i,k} = 1-a_{i,k}$ and $M$ is a sufficiently large number. 

The constraints \eqref{eqn:stateconstbin} on the binary variables $a_{i,k}$ ensure that the data is
compatible with at least one of the modes $i$ at each time $k$.
Note that in \eqref{eqn:stateconst0}-\eqref{eqn:stateconst1}, there exist bilinear terms in the form of
products of state and uncertainty variables. 
In what follows, we first apply a change of variables that renders these constraints linear. 
However, the constraints \eqref{eqn:uncertaintyconst3}-\eqref{eqn:uncertaintyconst4} are no longer linear in the new variables. 
Finally, we replace \eqref{eqn:uncertaintyconst3}-\eqref{eqn:uncertaintyconst4} 
with equivalent constraints that are linear in the new variables. 

Let us start by introducing the matrix variables $\mathbf{Z}_{i,k}^{\A}$ and $\mathbf{Z}_{i,k}^{\C}$ with entries:
\begin{equation} \label{eqn:chofv1}
\mathbf{Z}_{i,k}^{\tau}(m,l)  \doteq \hat{\pmb{\tau}}_i(m,l)\Delta_{i,k}^{\tau}(m,l)\x_k(l), \; \pmb{\tau} \in \{ \A,\C \}.
\end{equation} 
Then, we can rewrite constraints \eqref{eqn:stateconst0}-\eqref{eqn:stateconst1} as linear constraints in optimization variables:
\begin{subequations} \small
	\label{eqn:stateoutputconst}
	\begin{empheq}{align}
			& -\bar{a}_{i,k}M\mathbf{1} \leq \x_{k+1} - \A_{i}\x_k -\mathbf{Z}^{\A}_{i,k}\mathbf{1}- (\B_{i}+\B^\Delta_{i,k}) \mathbf{u}_k - \mathbf{f}_{i}- \f^\Delta_{i,k} \leq \bar{a}_{i,k}M\mathbf{1}, \label{eqn:state0} \\
			& -\bar{a}_{i,k}M{\mathbf{1}} \leq \mathbf{y}_k - \C_{i}\x_k-\mathbf{Z}^{\C}_{i,k}\mathbf{1} -  \ETA_k \leq \bar{a}_{i,k}M{\mathbf{1}}.  \label{eqn:output1}
	\end{empheq}
\end{subequations}
  Bound constraints for the new variables can be obtained using conic equivalences \cite{Bertsimas2006Tractable}.
In particular, constraints \eqref{eqn:uncertaintyconst} and \eqref{eqn:uncertaintyconst3}-\eqref{eqn:uncertaintyconst4} are equivalent to the following constraints in the new variables:
\begin{subequations}
 	\label{eqn:setmilpconst}
 	\begin{empheq}{align}
 		& \forall m,l \in \mathbb{Z}_n^+, \;  \forall o \in \mathbb{Z}_{n_y}^+:  \nonumber \\
 		& |\mathbf{Z}^{A}_{i,k}(m,l)| \leq |\hat{\A}_{i}(m,l)||\x_{k}(l)| \label{eqn:setmilpconst1} \\
 		& |\mathbf{Z}^{C}_{i,k}(o,l)| \leq |\hat{\C}_{i}(o,l)||\x_{
 			k}(l)|.\label{eqn:setmilpconst4}
 	\end{empheq}
\end{subequations}

\begin{proposition} \label{pro:mimilp}
	Given an SWA model $\G$ and an input-output trajectory $\{\uu_k,\yy_k\}_{k=0}^{N-1}$, assume there exists a large enough $M$ such that the following holds for all admissible states, uncertainty and noise values: 
	\begin{equation} \label{eqn:bigM1}
	\begin{aligned} 
	& | \x_{k+1} - (\mathbf{A}_{i}+\A^\Delta_{i,k}) \x_k -  (\B_{i}+\B^\Delta_{i,k}) \mathbf{u}_k- \mathbf{f}_{i}- \f^\Delta_{i,k}| \leq M\mathbf{1}, \\
	& |\mathbf{y}_k - (\C_{i}+\C^\Delta_{i,k})  \x_k - \ETA_k| \leq M{\mathbf{1}} ,
	\end{aligned}
	\end{equation}
	 then the model is invalidated if and only if the following problem is infeasible 
	\begin{align} \label{eqn:uncmilp} 
	\text{\emph{Find} } & \x_{k},\ETA_{k}, a_{i,k}, \mathbf{Z}^{\A}_{i,k},\Delta^{\B}_{i,k}, \mathbf{Z}^{\C}_{i,k}, \Delta^{\f}_{i,k}, \forall i\in\mathbb{Z}_s^+, \forall k\in \mathbb{Z}_{N-1}^0  \tag{P$_{MI}$}\\ 
	\text{\emph{s.t.} } & \left \{ \eqref{eqn:stateconstbin}\text{-}\eqref{eqn:uncertaintyconst2}, \eqref{eqn:stateoutputconst},\eqref{eqn:setmilpconst}
	 \right  \}, \; \forall k\in \mathbb{Z}_{N-1}^0,  \; \forall i\in \mathbb{Z}_s^+. \nonumber
	\end{align}
\end{proposition}
\begin{proof} If condition \eqref{eqn:bigM1} holds, then the satisfaction of state and output evolution constraints (cf. \eqref{eqn:SIMUSWA}) for mode $i$ at time $k$ is equivalent to $\bar{a}_{i,k}$ being zero in \eqref{eqn:stateconst0}-\eqref{eqn:stateconst1}. As the other constraints in the consistency set and \eqref{eqn:feasoptnonlineartot} are identical, the infeasibility of \eqref{eqn:feasoptnonlineartot} is equivalent to the emptiness of the consistency set, which by definition is equivalent to the invalidation of the model by data. Therefore it suffices to show the equivalence of  \eqref{eqn:feasoptnonlineartot} and \eqref{eqn:uncmilp}. Given a feasible point of \eqref{eqn:feasoptnonlineartot}, a feasible point of \eqref{eqn:uncmilp} can be constructed by applying the change of variables introduced by \eqref{eqn:chofv1}. For the other direction, let $\tilde{\x}_{k},\tilde{\ETA}_{k}, \tilde{a}_{i,k}, \tilde{\mathbf{Z}}^{\A}_{i,k}, \tilde{\Delta}^{\B}_{i,k}, \tilde{\mathbf{Z}}^{\C}_{i,k}, \tilde{\Delta}^{\f}_{i,k}$ be a feasible point of \eqref{eqn:uncmilp}. We will construct a feasible point $\check{\x}_k, \check{\ETA}_k, \check{a}_{i,k}, \check{\Delta}_{i,k}^{\A},  \check{\Delta}_{i,k}^{\B},  \check{\Delta}_{i,k}^{\C},  \check{\Delta}_{i,k}^{\f}$ for \eqref{eqn:feasoptnonlineartot}. Take $\check{a}_{i,k} = \tilde{a}_{i,k}$, which satisfies the constraints \eqref{eqn:stateconstbin} that is a common constraint. This means that at each time $k$, there is a unique mode $i^*_k$ such that $\check{a}_{i^*_k,k}=1$ and $\check{a}_{i,k}=0,\; \forall i \neq i^*_k$. Furthermore, take $\check{\x}_k = \tilde{\x}_k, \check{\ETA}_k = \tilde{\ETA}_k, \; \forall k \in \mathbb{Z}_{N-1}^0$, and $\check{\Delta}^{\B}_{i,k} = \tilde{\Delta}^{\B}_{i,k}, \check{\Delta}^{\f}_{i,k} = \tilde{\Delta}^{\f}_{i,k}, \forall i\in \mathbb{Z}_s^+, \forall k \in \mathbb{Z}_{N-1}^0$. Finally choose the elements of $\check{\Delta}_{i,k}^{\tau}, \tau \in \{\A,\C\}$ as follows:
\begin{equation} \label{eqn:Zs}
\check{\Delta}_{i,k}^{\tau}(m,l) = \begin{cases}
0 & \quad \text{if } \hat{\tau}_{i}(m,l)\check{\x}_k(l)=0\\
\frac{\tilde{\mathbf{Z}}^{\tau}_{i,k}\mathbf{1}(m)}{\hat{\tau}_i(m,l)\check{\x}_k(l)} & \quad \text{if } \hat{\tau}_{i}(m,l)\check{\x}_k(l)\neq0
\end{cases}
\end{equation}
Now, we check the feasibility of chosen variables by plugging them in \eqref{eqn:feasoptnonlineartot}.   

With the choice of $\check{\Delta}_{i,k}^{\tau}(m,l), \; \pmb{\tau} \in \{ \A, \C \}$ made in \eqref{eqn:Zs}, the following is true.
\begin{align*}
\tilde{\mathbf{Z}}^{\tau}_{i,k}\mathbf{1}(m) & = \sum_{l\in\mathbb{Z}^+_n}\tilde{\mathbf{Z}}^{\tau}_{i,k}(m,l) = \sum_{l\in\mathbb{Z}^+_n}\hat{\tau}_i(m,l)\check{\x}_k(l)\check{\Delta}_{i,k}^{\tau}(m,l)\\
& = [\pmb{\tau}_{i,k}^{\Delta}.\check{\x}_k](m) \implies \tilde{\mathbf{Z}}^{\tau}_{i,k}\mathbf{1} = \pmb{\tau}_{i,k}^{\Delta}.\check{\x}_k.
\end{align*}
This proves the equivalence of \eqref{eqn:stateoutputconst} and \eqref{eqn:stateconst0}-\eqref{eqn:stateconst1}. Constraints \eqref{eqn:stateconstbin}-\eqref{eqn:uncertaintyconst2} are trivially satisfied for our choice of feasible point. Now consider \eqref{eqn:setmilpconst1}, and plug in \eqref{eqn:Zs}, we have:
\begin{align*}
|\hat{\A}_i(m,l)\check{\x}_k(l)\check{\Delta}_{i,k}^{\A}(m,l)| = |\hat{\A}_i(m,l)||\check{\x}_k(l)||\check{\Delta}_{i,k}^{\A}(m,l)|\\ \leq |\hat{\A}_i(m,l)||\check{\x}_k(l)|
\implies |\check{\Delta}_{i,k}^{\A}(m,l)|\leq 1.
\end{align*}   
With the similar analysis for \eqref{eqn:setmilpconst4}, we can show that with our feasible point choice, from feasibility of \eqref{eqn:setmilpconst}, the feasibility of \eqref{eqn:uncertaintyconst3}-\eqref{eqn:uncertaintyconst4} is implied.
Therefore, the two problems are equivalent in the sense that the feasibility of one implies the feasibility of the other.
\end{proof}

We refer to an instance of the problem \eqref{eqn:uncmilp}
for a given SWA model $\G$ and an input-output trajectory $\{\uu_k,\yy_k\}_{k=0}^{N-1}$ 
as $\text{Feas}_{\G}(\{\uu_k,\yy_k\}_{k=0}^{N-1})$ in the
remainder of the paper. 
\begin{remark}
	Assumption \eqref{eqn:bigM1} for the existence of a large enough $M$ that bounds the state and output equations is trivially satisfied 
	when the admissible state set $\mathcal{X}$ is a bounded set. This assumption holds in practice, especially because physical systems often have bounded states. 
	If such $M$ exists, then it can be easily computed based on the admissible set of states, inputs and the model parameters. 
	On the other hand, if such an $M$ does not exist (admissible sets are unbounded), then this formulation does not exactly encode the model invalidation problem. 
	However, it is still possible to encode the model invalidation problem in an exact manner with a mixed integer linear programming problem by employing 
	a  convex hull reformulation, which is proposed in our earlier work \cite{Harirchi2015Model}. We resort to the big M formulation in this work because 
	for both model invalidation and $T$-detectability, its computational complexity is lower than other available techniques.
\end{remark}
Both the objective function and the constraints of the problem \eqref{eqn:uncmilp} are linear
except for the constraints in \eqref{eqn:setmilpconst}. However, it is straightforward to
transform these constraints into a set of linear constraints as shown in the following remark. 
\begin{remark}\label{rem:absmilp}
	Consider the first constraint in \eqref{eqn:setmilpconst}:
	This constraint is in the form of $|x|\leq C|y|, C \geq 0$ and both $x$ and $y$ are variables. Let us introduce a new variable $z=|y|$. Two cases can be considered: if $y\geq 0$, then $z=y$ and $z=-y$ vice versa. We also assign a binary variable $b$, which is 1, if $y\geq0$ and 0 otherwise. Then, $|x|\leq C|y|$ is enforced by the following linear constraints, where $M$ is an upper bound for $2|y|$:
	\begin{align*}
	-Cz \leq & x \leq Cz \\
	0 \leq & z-y \leq M (1 - b)\\
	0 \leq & z+y \leq Mb
	\end{align*}
\end{remark}
After this transformation, problem \eqref{eqn:uncmilp} becomes a MILP problem in feasibility form. 

\subsection{$T$-Detectability}
\begin{figure}[h!]
 	\centering
 	\includegraphics[width=5in]{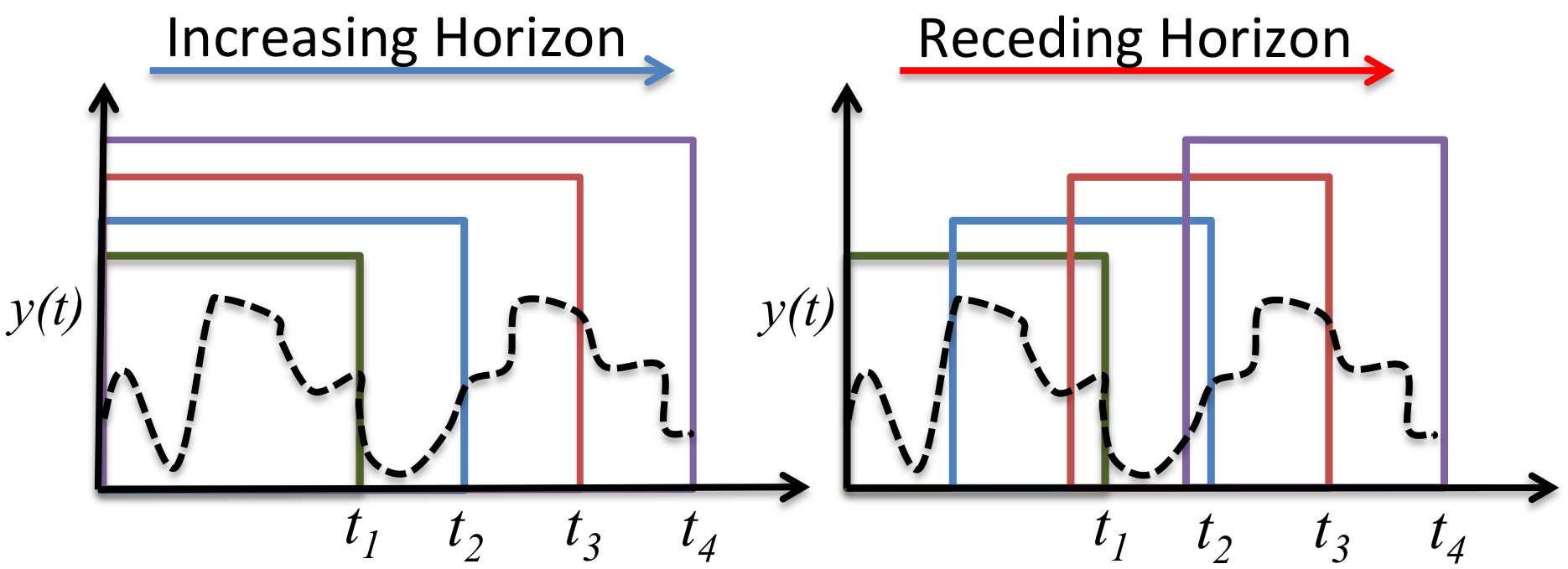}
 	\caption{Left: without $T$-detectability, the model invalidation has to be applied
 		on an increasing-size data. Right: with $T$-detectability, the model invalidation
 		is applied in a receding horizon manner.} \label{fig:recedhor}
 \end{figure}

In general, for detecting an arbitrary, unknown fault (or anomaly),
one needs to use all the data that is available, which leads to
larger and larger problems 
as time passes. On the other hand,
as alluded to earlier, the availability of a fault model can be exploited 
to develop more efficient fault detection schemes
while preserving the detection guarantees. 
We are interested in $T$-detectability of the pair of system and fault models as 
 a key property that enables us to apply the model invalidation approach for 
 fault detection in a receding horizon manner, rather than applying it on the entire time 
 horizon, while preserving all the guarantees for detection. Fig. \ref{fig:recedhor} 
 illustrates the motivation of introducing the concept of $T$-detectability.
 
  In this section, we develop a MILP characterization of 
 $T$-detectability property that can be verified off-line. Let us first make the following assumption on the fault models:
 \begin{assumption}\label{a:perm}
 	Faults are persistent, i.e., once occurred it will be permanent.
 \end{assumption}
In Section \ref{sec:gen}, we show how to relax Assumption \ref{a:perm}. 

The following proposition formalizes the fact that 
$T$-detectable faults can be detected with a receding horizon algorithm.

\begin{proposition}\label{prop:window} Given a $T$-detectable fault model $\G^f$ for an
	 SWA model $\G$, it is possible to detect the existence of a persistent fault evolving
	  through $\G^f$ by checking, at each time $k$, if $\text{Feas}_{\G}$ $(\{\uu_j,\yy_j\}_{j=k-T+1}^{k})$
	   is feasible or not.
\end{proposition}
\begin{proof} Let a fault occur at time $i^*$; that is, the input-output trajectory
	$ \{ \uu_j, \yy_j \}_{j\geq i^*}$ is generated by the fault model $\G^f$. 
	Because $\G^f$ is $T$-detectable, by Def. \ref{def:T_detec}, there exists
	$k^*\leq i^*+T-1$ such that $ \{ \uu_j, \yy_j \}_{j=i^*}^{k^*} \notin \mathcal{B}^{k^*-i^*+1}_{swa}(\G)$.
	By Proposition \ref{pro:mimilp}, this is equivalent to the existence of $k^*\leq i^*+T-1$ such that $\text{Feas}_{\G}(\{\uu_j,\yy_j\}_{j=i^*}^{k^*})$ is infeasible. Since 
	$[i^*,k^*] \subseteq  [k^*-T+1, k^*]$,  the infeasibility of 
	$\text{Feas}_{\G}(\{\uu_j,\yy_j\}_{j=i^*}^{k^*})$ implies the infeasibility of $\text{Feas}_{\G}(\{\uu_j,\yy_j\}_{j=k^*-T+1}^{k^*})$.
\end{proof}

Now that we have an efficient way to solve fault detection problem for a $T$-detectable fault,
the next question is: given an integer $T$, how to check whether a fault is $T$-detectable for a particular system or not? 
In what follows, we give a necessary and sufficient condition under which a fault $\G^f$ 
is $T$-detectable for a system $\G$ that can be verified by checking the feasibility of
MILP certificates.

In terms of consistency sets, $\G^f$ is $T$-detectable for $\G$ if for any $\{ \uu_k, \yy_k  \}_{k=0}^{T-1} \in \mathcal{B}^{T}_{swa}(\G^f)$, we 
have $\mathcal{T}^{*}_{\G}(\{ \uu_k, \yy_k  \}_{k=0}^{T-1}) = \emptyset$,
or equivalently for any $\{ \uu_k, \yy_k  \}_{k=0}^{T-1}$, the following holds:
\begin{equation} \label{eqn:consis}
\mathcal{T}^{*}_{\G}(\{ \uu_k, \yy_k  \}_{k=0}^{T-1}) \cap \mathcal{T}^{*}_{\G^f}(\{ \uu_k, \yy_k  \}_{k=0}^{T-1}) =  \emptyset. 
\end{equation}
Consider the following nonlinear feasibility problem:
\begin{subequations}
\begin{align} \label{eqn:TfeasoptMIP} \tag{19}
\text{Find } & \x_k, \bar{\x}_k,\ETA_k, \bar{\ETA}_k, \uu_k, \Delta_k, \bar{\Delta}_k, d_{i,j,k}, \forall
k \in \mathbb{Z}_{T-1}^0, \forall i \in \mathbb{Z}_{s}^+, \; \forall j \in \mathbb{Z}_{\bar{s}}^+ \qquad \qquad \qquad \qquad \nonumber \\
\text{s.t. } & \forall k \in \mathbb{Z}_{T-1}^0, \; \forall i \in \mathbb{Z}_{s}^+,\;  \forall j \in \mathbb{Z}_{\bar{s}}^+: \nonumber
\end{align} \vspace*{-0.3cm}
\begin{empheq}[left={ \empheqlbrack}]{align} 
& -\bar{d}_{i,j,k}M\mathbf{1} \leq \x_{k+1} - (\mathbf{A}_{i}+ \A^\Delta_{i,k}) \x_k - (\B_{i}+\B^\Delta_{i,k}) \mathbf{u}_k -\mathbf{f}_{i} - \f^\Delta_{i,k} \leq  \bar{d}_{i,j,k}M\mathbf{1}, \label{eqn:statecons1}\\
& -\bar{d}_{i,j,k}M\mathbf{1} \leq \bar{\x}_{k+1} - (\bar{\mathbf{A}}_{j}+\bar{\A}^{\Delta}_{j,k}) \bar{\x}_k-  (\bar{\B}_{j}+\bar{\B}^{\Delta}_{j,k})\uu_k -\bar{\mathbf{f}}_{j} -\bar{\f}^{\Delta}_{j,k} \leq \bar{d}_{i,j,k}M\mathbf{1},  \\
& -\bar{d}_{i,j,k}M\mathbf{1} \leq (\C_{i}+\C^\Delta_{i,k}) \x_k +\ETA_k  - (\bar{\C}_{j}+\bar{\C}^{\Delta}_{j,k})  \bar{\x}_k-\bar{\ETA}_k \leq \bar{d}_{i,j,k}M\mathbf{1}, \label{eqn:statecons2}\\
& \textstyle \sum_{i \in \mathbb{Z}_{s}^+} \sum_{j \in \mathbb{Z}_{\bar{s}}^+} d_{i,j,k} = 1, \; d_{i,j,k} \in \{0,1\}, \label{eqn:binconst}\\
& X_l \leq \x_k \leq X_u, \; \bar{X}_l \leq \bar{\x}_k \leq \bar{X}_u,  \\ 
& \epsilon_l \leq \ETA_k \leq \epsilon_u,\; \bar{\epsilon}_l \leq \bar{\ETA_k} \leq \bar{\epsilon}_u
,\; U_l \leq \uu_k \leq U_u, \label{eqn:bndconst}\\
& \text{for } \pmb{\tau}\in \{ \A, \C, \B,\f, \bar{\A}, \bar{\C}, \bar{\B},\bar{\f} \}:\\
& \pmb{\tau}^\Delta_{i,k} = \hat{\pmb{\tau}}_i \odot \Delta^{\tau}_{i,k}, \; |\Delta^{\tau}_{i,k}(m,l)| \leq 1,  \label{eqn:unc2} 
\end{empheq} 
\end{subequations}
where $M$ is a sufficiently large scalar and $\bar{d}_{i,j,k} = 1- d_{i,j,k}$. Clearly, the infeasibility of \eqref{eqn:TfeasoptMIP} is equivalent to \eqref{eqn:consis} being true, 
therefore, the $T$-detectability of $\G^f$ for $\G$. This is stated formally below.
\begin{proposition} \label{pro:Tdetect}
	Assume there exists a large enough $M$ such that the following holds for all admissible states, uncertainty and noise values for both system and fault models:
		\begin{equation} \label{eqn:bigM}
		\begin{aligned} 
		& | \x_{k+1} - (\mathbf{A}_{i}+\A^\Delta_{i,k}) \x_k -  (\B_{i}+\B^\Delta_{i,k}) \mathbf{u}_k- \mathbf{f}_{i}- \f^\Delta_{i,k}| \leq M\mathbf{1}, \\
		& | \bar{\x}_{k+1} - (\bar{\mathbf{A}}_{i}+\bar{\A}^\Delta_{i,k}) \bar{\x}_k -  (\bar{\B}_{i}+\bar{\B}^\Delta_{i,k}) \mathbf{u}_k- \bar{\mathbf{f}}_{i}- \bar{\f}^\Delta_{i,k}| \leq M\mathbf{1}, \\
		& |(\C_{i}+\C^\Delta_{i,k}) \x_k +\ETA_k - (\bar{\C}_{j}+\bar{\C}^{\Delta}_{j,k})  \bar{\x}_k-\bar{\ETA}_k| \leq M{\mathbf{1}} ,
		\end{aligned}
		\end{equation}
	then fault model $\G^f$ is $T$-detectable for the system model $\G$ if and only if problem
	 \eqref{eqn:TfeasoptMIP} is infeasible. 
\end{proposition}
\begin{proof}
The binary variables $d_{i,j,k}$ in \eqref{eqn:TfeasoptMIP} guarantee that at each time $k$,
there is a mode $j$ of the fault model $\G^f$ and 
a mode $i$ of the system model $\G$ such that the state evolution constraints for both models are satisfied and 
the outputs of the two models match. 
The rest of the proof is similar to the proof of Proposition \ref{pro:mimilp} and is omitted for brevity. 
\end{proof}

Following steps similar to those in the model invalidation subsection, we derive 
a MILP problem equivalent to Problem \eqref{eqn:TfeasoptMIP}.
In order to eliminate the bilinear terms in \eqref{eqn:TfeasoptMIP} that are shaped by the products of state and uncertainty variables,
we consider the change of variables introduced by \eqref{eqn:chofv1}. Similarly, we define ${\mathbf{Z}}^{\bar{A}}_{j,k}, {\mathbf{Z}}^{\bar{C}}_{j,k}$ for fault model as:
\begin{equation} \label{eqn:chofv2}
\mathbf{Z}_{j,k}^{\bar{\tau}}(m,l)  \doteq \hat{\bar{\pmb{\tau}}}_j(m,l)\Delta_{j,k}^{\bar\tau}(m,l)\bar{\x}_k(l), \; \bar{\pmb{\tau}} \in \{ \bar{\A},\bar{\C} \}.
\end{equation} 
Note that in contrast to model invalidation problem, inputs are variables for the $T$-detectability problem. Therefore, we need to introduce the following change of variables: 
\begin{equation} \label{eqn:chofv3}
\mathbf{Z}_{i,k}^{\tau}(m,p)  \doteq \hat{\pmb{\tau}}_i(m,p)\Delta_{i,k}^{\tau}(m,p)\uu_k(p), \; \pmb{\tau} \in \{ {\B},\bar{\B} \}.
\end{equation} 
where the indices cover all the entries of the corresponding matrices.
By plugging in the new variables we can equivalently state the constraints \eqref{eqn:statecons1}-\eqref{eqn:statecons2} by the following MILP constraints:
\begin{equation} \label{eqn:Tstateeqmilp} 
\begin{aligned}
  & -\bar{d}_{i,j,k}M\mathbf{1} \leq \x_{k+1} -\mathbf{A}_{i}\x_k- \mathbf{Z}^{A}_{i,k}\mathbf{1} - \B_{i}\mathbf{u}_k - \mathbf{Z}^{B}_{i,k}\mathbf{1}- \mathbf{f}_{i} - \f^\Delta_{i,k} \leq  \bar{d}_{i,j,k}M\mathbf{1}, \\
  & -\bar{d}_{i,j,k}M\mathbf{1} \leq \bar{\x}_{k+1} -\bar{\mathbf{A}}_{j}\bar{\x}_k- \mathbf{Z}^{\bar{A}}_{j,k}\mathbf{1} - \bar{\B}_{j}\mathbf{u}_k -  \mathbf{Z}^{\bar{B}}_{j,k}\mathbf{1}- \bar{\mathbf{f}}_{j} - \bar{\f}^\Delta_{j,k} \leq \bar{d}_{i,j,k}M\mathbf{1} .
\end{aligned}
\end{equation}

Now, we formulate constraint \eqref{eqn:unc2} in the new variables as MILP constraints, where $m,l$ encompass all the entries of the corresponding matrix.
\vspace{-0cm}
\begin{equation}\label{eqn:Tsetmilpp} 
\begin{aligned}
& |\mathbf{Z}^{\tau}_{i,k}(m,l)| \leq |\hat{\tau}_{i}(m,l)||\x_{k}(l)|, \tau \in \{\A, \C \}\\
& |\mathbf{Z}^{\bar{\tau}}_{j,k}(m,l)| \leq |\hat{\bar\tau}_{j}(m,l)||\bar{\x}_{k}(l)|, \bar{\tau} \in \{\bar{\A}, \bar{\C} \}  \\
& |\mathbf{Z}^{\tau}_{i,k}(m,l)| \leq |\hat{\tau}_{i}(m,l)||\uu_{
	k}(l)|, \tau \in \{\B, \bar{\B} \}\\
& \tau^\Delta_{i,k} = \hat{\tau}_i \odot \Delta^{\tau}_{i,k}, \; |\Delta^{\tau}_{i,k}(m,l)| \leq 1, \tau \in   \{\f, \bar{\f} \}.
\end{aligned}
\end{equation}
The absolute value constraints can be converted to 
 equivalent mixed integer linear constraints
as illustrated in Remark \ref{rem:absmilp}. 
\vspace{-0cm}
\begin{theorem} \label{thm:main}
	The fault model $\G^f$ is $T$-detectable for the system $\G$ if and only if 
	\eqref{eqn:Tmilp} is infeasible. 
\begin{equation} \label{eqn:Tmilp}  \tag{P$_T$}
\begin{aligned} 
\text{Find    } & \x_{k}, \bar{\x}_k, \uu_k, \ETA_{k}, \bar{\ETA}_k, d_{i,j,k}, \mathbf{Z}^{A}_{i,k}, \mathbf{Z}^{\bar{A}}_{j,k}, \mathbf{Z}^{B}_{i,k}, \mathbf{Z}^{\bar{B}}_{j,k}, \mathbf{Z}^{C}_{i,k}, \mathbf{Z}^{\bar{C}}_{j,k}, \Delta^{f}_{i,k}, \bar{\Delta}^{f}_{j,k} \\
\text{s. t.    } & \forall
k\in\mathbb{Z}_{T-1}^0,\; \forall i\in\mathbb{Z}_s^+, \; \forall j \in \mathbb{Z}_{\bar{s}}^+: \left \{ \eqref{eqn:binconst}-\eqref{eqn:bndconst},  \eqref{eqn:Tstateeqmilp}-\eqref{eqn:Tsetmilpp} \right  \} 
\end{aligned}	
\end{equation}
\end{theorem} 
The proof follows essentially the same reasoning as the proof of Proposition \ref{pro:mimilp}, i.e., by showing the feasibility of Problem \eqref{eqn:TfeasoptMIP} is equivalent to Problem \eqref{eqn:Tmilp}, and is omitted for brevity.


\subsection{Converse results for detectability}

For a given pair of fault and system models and a given integer $T$, Theorem \ref{thm:main} provides a necessary and sufficient condition
for verifying the $T$-detectability of the fault for the particular system. In general, to compute the 
minimum such $T$, one needs to start with $T=1$ and solve a sequence of problems
of the form \eqref{eqn:Tmilp} while incrementing $T$ until infeasibility is achieved.
Note that a finite $T$ does not necessarily exist for all types of faults. 
Showing that there exist no finite $T$ is equivalent to proving that the reachable sets
of two constrained switched systems always have a non-empty intersection. 
This is a notoriously hard problem given that many reachability problems involving 
constrained switched systems are not known to be decidable \cite{Egerstedt}.

In what follows we analyze $T$-detectability in a simplified setting  
with switched affine autonomous models $\G$, $\bar{\mathcal{G}}$
for the system and the fault, where we omit the inputs, uncertainty, noise
and admissible set constraints on the states. The analysis is based on establishing a
relation between $T$-detectability and observability of switched systems. 
We first recall some results on observability of switched systems from \cite{Babaali2003Pathwise} . 

Consider an autonomous switched system, $\G$ with the following dynamics:
\begin{equation} \label{eqn:simple_sw}
\begin{aligned} 
\x_{k+1}  	& = \A_{\sigma_k} \x_k \\ 
\mathbf{y}_k  & = \C_{\sigma_k}  \x_k
\end{aligned}, \sigma_k \in \mathbb{Z}_s^+,
\end{equation}
where the state $\x \in \mathbb{R}^n$ and the output $\mathbf{y}\in \mathbb{R}$.
Given a mode sequence $\sigma_{0:N}$, the observability matrix
corresponding to this mode sequence is defined as
\begin{equation}
\mathcal{O}_{\G}(\sigma_{0:N}) \doteq \begin{bmatrix}
\C_{\sigma_0} \\
\C_{\sigma_1}\A_{\sigma_0} \\
\vdots \\
\C_{\sigma_N}\A_{\sigma_{N-1}}\hdots\A_{\sigma_{0}} 
\end{bmatrix}.
\end{equation}
The following result is from \cite{Babaali2003Pathwise} (cf. Theorem 1 and Lemma 6 in \cite{Babaali2003Pathwise}).

\begin{lemma} \label{lem:baba} There exists a large enough finite number $\mathcal{N} \doteq \mathcal{N}(s,n)$ that depends on the number of modes $s$ and the state dimension $n$ such that if for a mode sequence $\sigma_{0:\mathcal{N}-1}$ of length $\mathcal{N}$, we have $\text{rank}(\mathcal{O}_{\G}(\sigma_{0:\mathcal{N}-1})) <n$, then for any arbitrary length there exists a mode sequence $\hat{\sigma}$ of that length such that:
	\begin{equation}
	 \text{{\em Range}}(\mathcal{O}_{\G}(\hat{\sigma})) \subseteq  \text{{\em Range}}(\mathcal{O}_{\G}(\sigma_{0:\mathcal{N}-1})).
	\end{equation}
\end{lemma}
The number $\mathcal{N}(s,n)$ can be recursively computed given the values of $s$ and $n$. We refer the reader to \cite{Babaali2003Pathwise} for its computation.

Now consider a pair of autonomous switched affine models $\G$, $\bar{\mathcal{G}}$ of the form:
\begin{equation} \label{eqn:simpswa}
\G: \begin{cases} \x_{k+1}  	 = \A_{\sigma_k} \x_k + \mathbf{f}  \\ 
\mathbf{y}_k   = \C_{\sigma_k}  \x_k\\
\sigma_k \in \mathbb{Z}_s^+
\end{cases},
\bar{\mathcal{G}}: \begin{cases} \bar{\x}_{k+1}  = \bar{\A}_{\bar{\sigma}_k} \bar{\x}_k + \bar{\mathbf{f}}  \\ 
\bar{\mathbf{y}}_k  = \bar{\C}_{\bar{\sigma}_k}  \bar{\x}_k\\
\bar{\sigma}_k \in \mathbb{Z}_{\bar{s}}^+ 
\end{cases}
\end{equation}

Our goal is to provide necessary and sufficient conditions for the existence of a finite $T$ for $T$-detectability of $\G$ and $\bar{\mathcal{G}}$. First we define equivalent linear models both for the system and the fault by augmenting their states with affine terms. These equivalent systems, denoted by $\G_{aug}$ and $\bar{\mathcal{G}}_{aug}$, have the following form:
\begin{equation}
\begin{aligned} 
\begin{bmatrix}
\x_{k+1}\\
\f
\end{bmatrix}  	& = \begin{bmatrix}
\A_{\sigma_k} & \mathbf{I}\\
\mathbf{0} & \mathbf{I}
\end{bmatrix} \begin{bmatrix}
\x_{k}\\
\f
\end{bmatrix} \\ 
\mathbf{y}_k  & = \begin{bmatrix}
\C_{\sigma_k} & \mathbf{0}
\end{bmatrix} \begin{bmatrix}
\x_{k}\\
\f
\end{bmatrix}
\end{aligned}, \sigma_k \in \mathbb{Z}_s^+.
\end{equation}
In order to check if the two models $\G$, $\bar{\mathcal{G}}$ of the form \eqref{eqn:simpswa} are not $T$-detectable, it suffices to check whether or not there exist initial conditions $\x_0, \bar \x_0$ and mode sequences $\sigma, \bar{\sigma}$ of length $T$ such that the following is true:
\begin{equation} \label{eqn:notdetect}
\mathcal{O}_{\G_{aug}}(\sigma) \begin{bmatrix}
\x_0 \\
\f
\end{bmatrix} - \mathcal{O}_{\bar{\mathcal{G}}_{aug}}(\bar \sigma) \begin{bmatrix}
\bar \x_0 \\
\bar \f
\end{bmatrix} = \mathbf{0},
\end{equation}
where this equation is obtained by setting the outputs of the system and fault models to be equal. Therefore, the set of initial conditions $\x_0, \bar \x_0$ and mode sequences $\sigma, \bar{\sigma}$ of length $T$ satisfying \eqref{eqn:notdetect} are essentially a projection of the feasible set of the problem \eqref{eqn:Tmilp} onto these variables.

Let us define a concatenated switched system:
\begin{equation}\label{eqn:concatenated}
\tilde{\mathcal{G}}: \begin{cases}
\begin{bmatrix}
\x_{k+1}\\
\f\\
\bar{\x}_{k+1}\\
\bar{\f}
\end{bmatrix}  	= \begin{bmatrix}
\A_{\sigma_k} & \mathbf{I} & \mathbf{0} & \mathbf{0} \\
\mathbf{0} & \mathbf{I} & \mathbf{0} & \mathbf{0}\\
\mathbf{0} & \mathbf{0} & \bar{\A}_{\bar{\sigma}_k}& \mathbf{I}\\
\mathbf{0} & \mathbf{0} &  \mathbf{0} & \mathbf{I}
\end{bmatrix} \begin{bmatrix}
\x_{k}\\
\f\\
\bar{\x}_{k}\\
\bar{\f}
\end{bmatrix}  	 \\ 
\mathbf{y}_k  = \begin{bmatrix}
\C_{\sigma_k} & \mathbf{0} & -\bar{\C}_{\bar{\sigma}_k} &  \mathbf{0}
\end{bmatrix} \begin{bmatrix}
\x_{k}\\
\f\\
\bar{\x}_{k}\\
\bar{\f}
\end{bmatrix},
\begin{array}{l} 
\sigma_k \in \mathbb{Z}_s^+\\ 
\bar{\sigma}_k \in \mathbb{Z}_{\bar s}^+ 
\end{array}
\end{cases}
\end{equation}
Also define $\tilde{\sigma}$ as a path of system $\tilde{\mathcal{G}}$ that is shaped by choosing paths $\sigma, \bar{\sigma}$ for $\G_{aug} , \bar{\mathcal{G}}_{aug}$, respectively. In other words, $\tilde{\sigma}$ take values from $\tilde{\Sigma} = \{(i,j)\mid i\in \mathbb{Z}_s, j\in \mathbb{Z}_{\bar{s}} \}$, with $|\tilde{\Sigma}| =s\bar{s} $. We can write the observability matrix of the concatenated system for switching signal $\tilde{\sigma}$ as follows:
\begin{equation}
\mathcal{O}_{\tilde{\mathcal{G}}}(\tilde{\sigma}) = \begin{bmatrix}
\mathcal{O}_{\G_{aug}}(\sigma)  & -\mathcal{O}_{\bar{\mathcal{G}}_{aug}}(\bar{\sigma}) 
\end{bmatrix}.
\end{equation}
Then, we can write the condition for not being $T$-detectable (i.e. \eqref{eqn:notdetect}) as the existence of $\x_0$, $\bar{\x}_0$ and $\tilde{\sigma}$ such that:
\begin{equation} \label{eqn:notdetect1}
\mathcal{O}_{\tilde{\mathcal{G}}}(\tilde{\sigma})\begin{bmatrix}
\x_0^\intercal &  \f^\intercal & \bar{\x}_0^\intercal & \bar {\f}^\intercal
\end{bmatrix}^\intercal = \mathbf{0}.
\end{equation}

\begin{theorem} \label{prop:NSaffine} Two switched affine models $\G$, $\bar{\mathcal{G}}$ of the form \eqref{eqn:simpswa} are {\bf not} $T$-detectable for any finite $T$ if and only if 
the problem \eqref{eqn:Tmilp} is feasible for $T=\mathcal{N}(s\bar{s}, 4n)$.
\end{theorem}
\vspace{-0.4cm}
\begin{proof} If \eqref{eqn:Tmilp} is infeasible for $T=\mathcal{N}(s\bar{s}, 4n)$, then the pair is trivially $T$-detectable. If \eqref{eqn:Tmilp} is feasible for $T=\mathcal{N}(s\bar{s}, 4n)$, then there exists initial conditions $\x_0, \bar \x_0$ and mode sequences $\sigma, \bar \sigma$ of length $\mathcal{N}(s\bar{s}, 4n)$ such that \eqref{eqn:notdetect1} holds. By Lemma \ref{lem:baba}, for the switched system \eqref{eqn:concatenated}, if $\text{rank}(\mathcal{O}_{\tilde{G}}(\sigma_1))<4n$ for a $\tilde{\sigma}_1$ of length $\mathcal{N}(s\bar{s}, 4n)$, then for any arbitrary time horizon $T$, there exists a mode sequence $\tilde{\sigma}_2$ of length $T$ such that $\text{Range}(\mathcal{O}_{\tilde{\mathcal{G}}}(\tilde{\sigma}_2)) \subseteq \text{Range}(\mathcal{O}_{\tilde{\mathcal{G}}}(\tilde{\sigma}_1))$, which implies that $\text{Null}(\mathcal{O}_{\tilde{\mathcal{G}}}(\tilde{\sigma}_1)) \subseteq \text{Null}(\mathcal{O}_{\tilde{\mathcal{G}}}(\tilde{\sigma_2}))$. Therefore if there exist initial conditions $\x_0$, $\bar{\x}_0$ that satisfy 
	\begin{equation}
	\begin{bmatrix}
	\x_0^\intercal &  \f^\intercal & \bar{\x_0}^\intercal & \bar {\f}^\intercal
	\end{bmatrix}^\intercal \in \text{Null}(\mathcal{O}_{\tilde{\mathcal{G}}}(\tilde{\sigma}_1)),
	\end{equation}
we have
	\begin{equation}
	\begin{bmatrix}
	\x_0^\intercal &  \f^\intercal & \bar{\x_0}^\intercal & \bar {\f}^\intercal
	\end{bmatrix}^\intercal \in \text{Null}(\mathcal{O}_{\tilde{\mathcal{G}}}(\tilde{\sigma}_2)).
	\end{equation}
Since $\tilde{\sigma}_2$ is of arbitrary length, and by noting that any mode sequence for \eqref{eqn:concatenated} can be mapped to mode sequences for $\G$ and $\bar{\mathcal{G}}$ of the same length, we conclude that the two models $\G$ and $\bar{\mathcal{G}}$ are not $T$-detectable for any finite $T$. 
\end{proof}

\begin{coro} Two affine models $G$, $\bar{G}$ of the form \eqref{eq:affine_mod} are {\bf not} $T$-detectable for any finite $T$ 
if and only if the problem \eqref{eqn:Tmilp} is feasible for $T=2n+1$.
\end{coro}
\vspace{-0.4cm}
\begin{proof} The feasible set of the problem \eqref{eqn:Tmilp} in this case projected onto the initial states is the set of $\x_0$, $\bar\x_0$ satisfying: 
\begin{align} 
\nonumber \begin{bmatrix}
\C  & - \bar{\C} 
\end{bmatrix}\begin{bmatrix} \x_0 \\ \bar{\x}_0\end{bmatrix} &= \mathbf{0} \\
\nonumber \begin{bmatrix} \mathcal{O}_G(2n) & -\mathcal{O}_{\bar{G}}(2n)\end{bmatrix} \left ( \begin{bmatrix}\A -I & 0 \\ 0 & \bar\A-I \end{bmatrix} \begin{bmatrix}  \x_0 \\ \bar\x_0 \end{bmatrix} + \begin{bmatrix}
\f \\
\bar{\f}
\end{bmatrix} \right ) &=\mathbf{0},
\end{align}
where $\mathcal{O}_G(2n)$ denotes the observability matrix associated with the $(\A, \C)$ pair for $2n$ steps.  
The result now follows by noting that the maximum rank of the matrix $[\mathcal{O}_G(2n) \; -\mathcal{O}_{\bar{G}}(2n)]$ is attained at $2n$.
\end{proof}


\subsection{Weak Detectability}
$T$-detectability, although useful to reduce the complexity of the proposed fault detection scheme,
can be rather strong in the sense that it is not hard to encounter faults in real applications
that are not $T$-detectable for any finite $T$. In this section, inspired by the indicators in
discrete-event systems diagnosis problems \cite{Sampath1995Diagnosability}, we propose
the notion of weak detectability that incorporates language constraints on the hidden mode
in the fault model. Since these indicators further restrict the behavior of the fault, it enlarges
the class of faults that can be detected within a finite horizon.

Let us motivate weak detectability and indicators with an example.
Consider a building radiant system with two modes that represent 
a controlled valve being open (mode 1) or closed (mode 2).
When the control logic is not modeled, this model allows 
arbitrary switching between the two modes. Now consider
a problem with the valve, which prevents it from being totally
closed. Therefore the fault model consists of two modes:
open and half-open. This fault model is not $T$-detectable
because if the valve is never commanded to close (i.e., 
system always operates in mode 1), it is not possible to
differentiate healthy system from the faulty one. 
On the other hand, this fault only affects the operation
of the system when the valve is commanded close so it 
will not be of interest to detect it if the valve is never
required to be closed. That is, the fault is relevant only if
mode 2 becomes active. In order to incorporate
such information in the fault models, we introduce 
indicators that shrink the behavior set of the SWA fault models
by restricting their
allowable mode sequences. 
Let us now formally define what we mean by an indicator.

\begin{definition} \label{def:indicators}
	(indicator) Given a fault model $\G^f$ with $s$ modes, 
	let $(\mathbb{Z}_s^+)^W$ be the set of all length-$W$ mode sequences.
	 An \emph{indicator} 
	   $\mathcal{I}$ is a subset of $(\mathbb{Z}_s^+)^W$ for some $W$. 
	   The \emph{fault-indicator} model $\G^f_{\mathcal{I}}$ is the fault model  
	   $\G^f$ whose mode sequences are restricted to 
	   $\mathcal{I}$ on the first $W$ time steps after the fault occurs.
\end{definition}

Indicators can be compactly defined using bounded regular languages or
fixed-length prefixes of some linear temporal logic formula \cite{Baier2008Principles}.
We also introduce a special class of indicators 
that we represent by tuples
\begin{align*}
	\mathcal{I}_t \doteq (\mathcal{S},W,m,O), O \in \{>,=,<\},
\end{align*}
where $\mathcal{I}_t$ consists of mode sequences of length $W$  
that contain modes from the set $\mathcal{S} \subseteq \mathbb{Z}_s^+$ 
more than ($O$ is $>$), exactly ($O$ is $=$), less than ($O$ is $<$) $m$ times. 
We have found this representation to be convenient in various 
applications we considered.

The behaviors of fault-indicator models are defined similarly.
\begin{definition}\label{def:f_ind_b}
	(fault-indicator behavior) The length-$N$ behavior of the 
	fault-indicator model $\G^f_{\mathcal{I}}$ is defined as:
	\begin{align*}
	&\mathcal{B}^N_{FI}(\G^f_{\mathcal{I}})  \doteq \big \{ \{\uu_k, \yy_k\}_{k=0}^{N-1} \mid 
	\uu_k \in \bar{\mathcal{U}}  \text{ and } \exists \x_k\in \bar{\mathcal{X}},  \ETA_k \in \bar{\mathcal{E}}, \Delta_k \in \bar{\Omega}, \{\sigma_k\}_{k=0}^{W-1}\in 
	\mathcal{I}, \text{ s.t. } \eqref{eqn:SIMUSWA} \text{ holds } \\
	& \text{  for fault model parameters} \big \}.
	\end{align*}
\end{definition}
In Def. \ref{def:f_ind_b}, we assume that $N \geq W$. It is possible to relax this assumption by considering
length-$N$ prefixes of sequences in $\mathcal{I}$. Now we are ready to define $T$-weak detectability.
\begin{definition} \label{def:weak}
	($T$-weak detectability) The fault model $\G^f$ is {\em $T$-weak detectable} with indicator $\mathcal{I}$
	 for the system model $\G$, if the fault-indicator model $\G_{\mathcal{I}}^f$ is $T$-detectable
	  for $\G$, that is, the following holds:
	\begin{equation}
	\mathcal{B}^T_{swa}(\G) \cap \mathcal{B}^T_{FI}(\G^f_{\mathcal{I}}) = \emptyset.
	\end{equation}  
\end{definition}
 Checking $T$-weak detectability also reduces to an equivalent MILP problem. 
 Recall that the binary variable $d_{i,j,k}$ in \eqref{eqn:Tmilp} indicates if mode $i$ of the system and
 mode $j$ of the fault are active at time $k$, 
 hence $\textstyle \sum_{i\in \mathbb{Z}_s^+} d_{i,j,k} = 1$ indicates that mode $j$ of the fault
  is active at time $k$.
 Then, if the mode sequence of the fault is $w_m \doteq \{\sigma^m_0 \hdots \sigma^m_{W-1} \} \in \mathcal{I}$,
 the following constraint needs to be enforced:
\begin{equation} \label{eqn:indicatormip}
\textstyle \sum_{i\in \mathbb{Z}^+_s}(d_{i,\sigma^m_0,0} + d_{i,\sigma^m_1,1} + \hdots +
 d_{i,\sigma^m_{W-1},W-1}) = W.
\end{equation}
Since the mode sequence of the fault should follow at least one $w_m \in \mathcal{I}$,
we have
 \begin{equation}\label{eqn:MIPconst}
b_m\big (\textstyle\sum_{i\in \mathbb{Z}^+_s} (d_{i,\sigma^m_0,0} + \hdots + d_{i,\sigma^m_{W-1},W-1}) -W \big ) = 0
\end{equation}
with $\textstyle \sum_{m\in \mathbb{Z}^+_\ell} b_m \geq 1$, where $b_m$ are new binary variables and $\ell \doteq | \mathcal{I} |$.
Finally, the change of variables, $d^m_{i,\sigma^m_k,k} \doteq b_m d_{i, \sigma^m_k,k}$, and extra constraints, $d^m_{i, \sigma^m_k,k} \leq b_m$, converts \eqref{eqn:MIPconst} to the equivalent MILP constraints :
\begin{equation} \label{eqn:MILPIND}
\begin{aligned}
& \textstyle\sum_{i\in \mathbb{Z}_s^+} d^m_{i,\sigma^m_k,k} -b_m W =0,\\
& d^m_{i, \sigma^m_k,k} \leq b_m, \; \textstyle \sum_{m\in \mathbb{Z}^+_\ell} b_m \geq 1.
\end{aligned}
\end{equation}
In general, structured indicators can be expressed by much less number of MILP constraints.
For instance the indicators of form $\mathcal{I}_t = (\mathcal{S},W,m,>)$ can be expressed by a single MILP
constraint:
\begin{equation}
\textstyle \sum_{i\in \mathbb{Z}^+_s} \sum_{j\in \mathcal{S}} \sum_{k\in \mathbb{Z}^+_W}  d_{i,j,k} \geq m.
\end{equation}

\subsection{Overall Fault Detection Scheme}
The overall fault detection scheme consists of an offline step and an online step. 
Given a system and a fault model (possibly with indicators), we first analyze $T$-(weak) detectability
running Algorithm \ref{Alg:findT} until it returns a $T$ or maximum number of iterations are reached.
\begin{algorithm}[h!]
	\caption{Calculating $T$} \label{Alg:findT}
	Input: $T_0, \G, \G^f_\mathcal{I}$
	\begin{enumerate}
		\item Initialize $T = T_0$. 
		\item Set $\mathcal{I}'$ to length-$T$ prefixes of $\mathcal{I}$
		\item Solve feasibility problem \eqref{eqn:Tmilp} with constraints \eqref{eqn:MILPIND} for $\mathcal{I}'$
		\item If  \eqref{eqn:Tmilp} is feasible:
		\begin{itemize}
			\item $T = T+1$ and go to step 2
		\end{itemize}
		Else
		\begin{itemize}
			\item Return $T$
		\end{itemize}   	
	\end{enumerate}
\end{algorithm}

Assuming that Algorithm \ref{Alg:findT} returns a finite $T$, at run-time,
we solve the problem \eqref{eqn:uncmilp} using data from the last $T$ steps at each time.
By construction, any occurrence of the fault is guaranteed to be detected
since it renders \eqref{eqn:uncmilp} infeasible and 
a false alarm never occurs. In case, $T$-detectability cannot be verified for a 
finite $T$, one can still use \eqref{eqn:uncmilp} in a receding horizon manner
but the occurrence of the fault can be missed in this case unless one
uses increasing amounts of data at each step. It is also worth noting that
if an unmodeled fault occurs, it can still be detected in this manner; i.e.,
any detection corresponds to an anomaly in the behavior of the system.
\vspace{-0.15cm}
\section{Generalizations}\label{sec:gen}
\vspace{-0.15cm}
Thanks to the generality of the modeling framework, the fault detection scheme 
presented in this paper 
can be applied to several other problems without any modifications. 
We briefly mention some of these problems in this section.
\vspace{-0.15cm}
\subsection{Attack Detection}
\vspace{-0.15cm}
Recently, there has been a considerable interest in attack detection \cite{pasqualetti2013attack} and 
attack resilient estimation \cite{Shoukry2014Secure,Chong2015Observability} for cyber-physical systems. 
One can argue that fault detection and attack detection are different problems
since, in general, an attacker, in addition to harming the system, has the intention of
not being noticed whereas the latter is not a concern for faults. 
However, since the proposed framework takes a worst-case approach, 
it is also well suited for attack detection. The only slight difference is that
attack models are typically subject to more uncertainty than fault models. 

Consider the following linear state-space model subject to sensor attack taken from \cite{Shoukry2014Secure}:
\begin{equation} \label{eqn:attackmdl}
\begin{aligned}
\x_{k+1} &= \A\x_k+\B\uu_k \\
\yy_k &= \C \x_k + \mathbf{a}_k + \ETA_k 
\end{aligned},
\end{equation}
where $\mathbf{a}_k \in \mathbb{R}^{n_y}$ is the attack vector at time instance $k$. 

The system model corresponds to the case where $\mathbf{a}_k=0$ in 
\eqref{eqn:attackmdl}. Take as the fault (i.e., attack) model as a
switched linear system with $\bar{s} = \textstyle \sum_{i\in \mathbb{Z}_a^+} {{n_y}\choose{i}}$ modes,
where each mode represents one possible attack scenario 
in which at most $a$ sensors are attacked. If $\mathcal{A} \subseteq \mathbb{Z}^+_{n_y}$, 
$|\mathcal{A}|\leq a$, is the set of sensors attacked in mode $j$, then
mode $j$ has the form \eqref{eqn:attackmdl}, with uncertain parameters $\mathbf{a}_k^i \neq 0$ 
(or, $|\mathbf{a}_k^i| \geq \epsilon$ if a minimum attack strength $\epsilon$ is known), 
$i\in \mathcal{A}$. Conditions like ``the same set of sensors are attacked persistently"
can be included using indicators that restrict the switching between attack modes. 
If this attack model can be shown to be $T$-detectable, running model invalidation
problem in a receding horizon manner detects the attacks. Moreover, if additional
observability conditions hold such as $s$-sparse observability in \cite{Shoukry2014Secure} with $s=2a$,
model invalidation problem can be used for resilient state estimation 
with an additional constraint that limits the number of attacked sensors.
One can also precisely estimate the trade-off between the noise level
and attack magnitude by solving $T$-detectability problems with appropriately defined objectives.
Note that our framework trivially extends to the case where the system itself
is switched or there are other types of attacks if the goal is detecting the attacks. 
Obtaining conditions for state estimation
guarantees in this general setting is left for future work.
\vspace{-0.15cm}
\subsection{Cascaded Faults}
\vspace{-0.15cm}
By using indicators and the fact that the proposed framework
handles switched systems, it is possible to relax  Assumption \ref{a:perm} on the persistency of the faults. 
In general, a system is subject to multiple faults and these
faults can occur in a cascaded fashion. Assume each such fault
is represented by an SWA model. It is possible to define a new SWA model
whose number of modes is the sum of the number of modes of the system model
and each fault model. Then, the cascading faults can be represented by restricting
the mode sequences of this new SWA model with indicators that
capture potential a priori knowledge on the order of failures. 
Therefore, one can apply the proposed framework to analyze 
$T$-detectability of multiple
cascading faults and to detect non-permanent
faults.

\vspace{-0.15cm}
\section{Illustrative Examples} \label{sec:example}
\vspace{-0.15cm}
In this section, we consider one set of numerical examples and an example motivated by building radiant 
systems to illustrate the efficacy of the methods proposed in this paper. All the examples are implemented 
on a 3.5 GHz machine with 32 GB of memory running Ubuntu. For the implementation of model invalidation 
approach and finding $T$ for $T-$detectability, we utilized Yalmip \cite{YALMIP} and CPLEX \cite{cplex}. All 
the approaches and examples are implemented in Matlab, and are available with 
MI4Hybrid\footnote{\url{https://github.com/data-dynamics/MI4Hybrid}} toolbox. 

\subsection{Numerical Results}
\begin{figure}[h!] 
	\centering
	\includegraphics[width=5in]{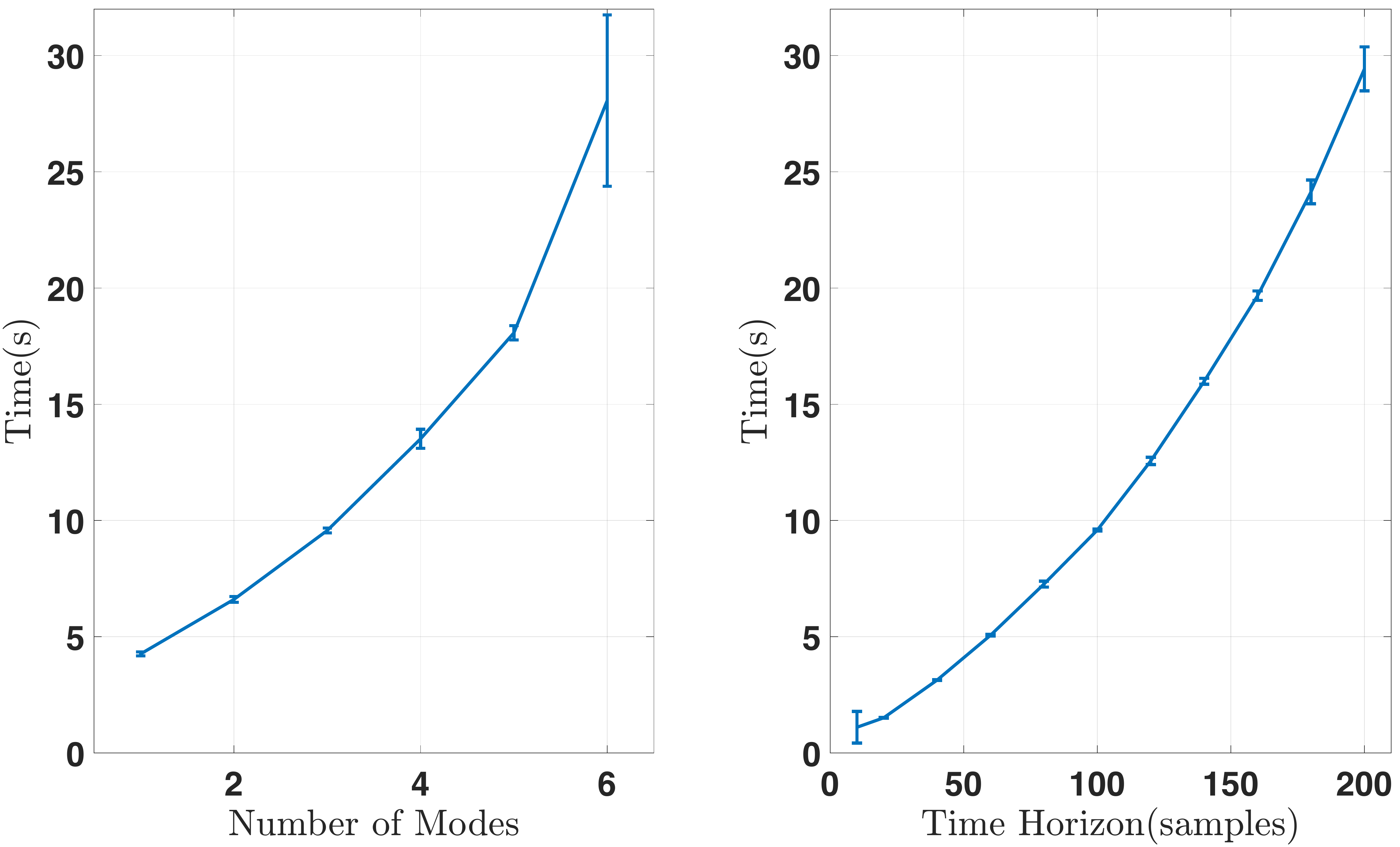}
	\caption{The average execution time for invalidation of data generated by $\G^f$ for $\G_3$ on various 
		time horizons (left), and for $\G_t$ on a fixed time horizon of length 100 for different values of $t$ (right).} 
	\label{fig:run_time}
\end{figure}

In this section, we illustrate the efficiency of the model invalidation approach proposed in this paper in terms 
of execution-time for various time-horizons and number of modes of the system. Consider a hidden-mode switched affine model, 
$\G_6$, with admissible sets $\X = \{ \x \mid\| \x \| \leq 11 \}$, $\U = \{ \uu \mid\| 
\uu \| \leq 1000 \}$ and $\E = \{ \ETA \mid\| \ETA \| \leq 0.1 \}$. We assume there is no process noise and 
no uncertainty in the parameters of the system. $\G_6$ has three states and six modes. We assume a fixed 
$B =[1 \;\; 0 \;\; 1]^T$ and $C = [1 \;\; 1 \;\; 1]$ for all modes. The system matrices of the modes are:

\begin{align*} 
& \A_1 = \begin{bmatrix}
0.5 & 0.5 & 0.5\\
 0.1 & -0.2 & 0.5\\
-0.4 & 0.6 & 0.2
\end{bmatrix}, \f_1 = \begin{bmatrix}
1 \\ 0 \\ 0
\end{bmatrix},  \A_2 = \begin{bmatrix}
0.5 & 0.5 & 0.5\\
-0.3 & -0.2 & 0.3 \\
0.1 & -0.3 & -0.5
\end{bmatrix},  \f_2 = \begin{bmatrix}
0 \\ 1 \\ 0
\end{bmatrix} \nonumber \\
& \A_3 = \begin{bmatrix}
0.5 & 0.2 & 0.6\\
0.2 & -0.2 & 0.2\\
-0.9 & 0.7 & 0.1
\end{bmatrix},  \f_3 = \begin{bmatrix}
0 \\ 0 \\ 1
\end{bmatrix},  \A_4 = \begin{bmatrix}
-0.5 & 0.5 & 0.8\\
0.1 & -0.2 & -0.6\\
0.2 & -0.6 & 0.3
\end{bmatrix}, \f_4 = \begin{bmatrix}
1 \\ 1 \\ 0
\end{bmatrix} \nonumber \\
& \A_5 = \begin{bmatrix}
0.8 & 0.5 & 0.2\\
-0.1 & 0.2 & -0.3\\
0.5 & 0.4 & -0.1
\end{bmatrix}, \f_5 = \begin{bmatrix}
0 \\ 1 \\ 1
\end{bmatrix},  \A_6 = \begin{bmatrix}
-0.3 & 0.8 & -0.1\\
0.4 & -0.1 & 0.3\\
0.9 & -0.2 & 0.6
\end{bmatrix}, \f_6 = \begin{bmatrix}
1 \\ 0 \\ 1
\end{bmatrix}. 
\end{align*} 
Let us define, $\G_t, t = 1, \hdots ,6$, to be a system with the first $t$ modes of $\G_6$. In addition, 
suppose the fault, $\G^f$, is an affine model with the following system matrices:
\begin{align} \label{eqn:numexfault}
& \A^f = \begin{bmatrix}
0.8 & 0.7 & 0.6\\
0.1 & -0.2 & 0.3\\
-0.4 & 0.3 & -0.2
\end{bmatrix}, \B^f = \begin{bmatrix}
1 \\ 0 \\ 0
\end{bmatrix}, \f^f = \begin{bmatrix}
2.5 \\ 2.5 \\ 2.5
\end{bmatrix}.
\end{align}
We first generate input-output trajectories from $\G^f$ of various lengths. 
We then fix the length at 100 samples and investigate the effect of increasing 
the number of modes of the system. 
Both examples are repeated for 10 times for randomly 
generated input and noise sequences. The results for applying model invalidation
in this setting are illustrated in Fig. \ref{fig:run_time}. Although the execution-time
scales relatively reasonably with increasing horizon and number of modes,
depending on the time-scale of the system dynamics,
solving only short horizon problems might be feasible
for real-time implementations, motivating $T$-detectability.

\begin{figure}[h!] 
	\centering
	\includegraphics[width=5in]{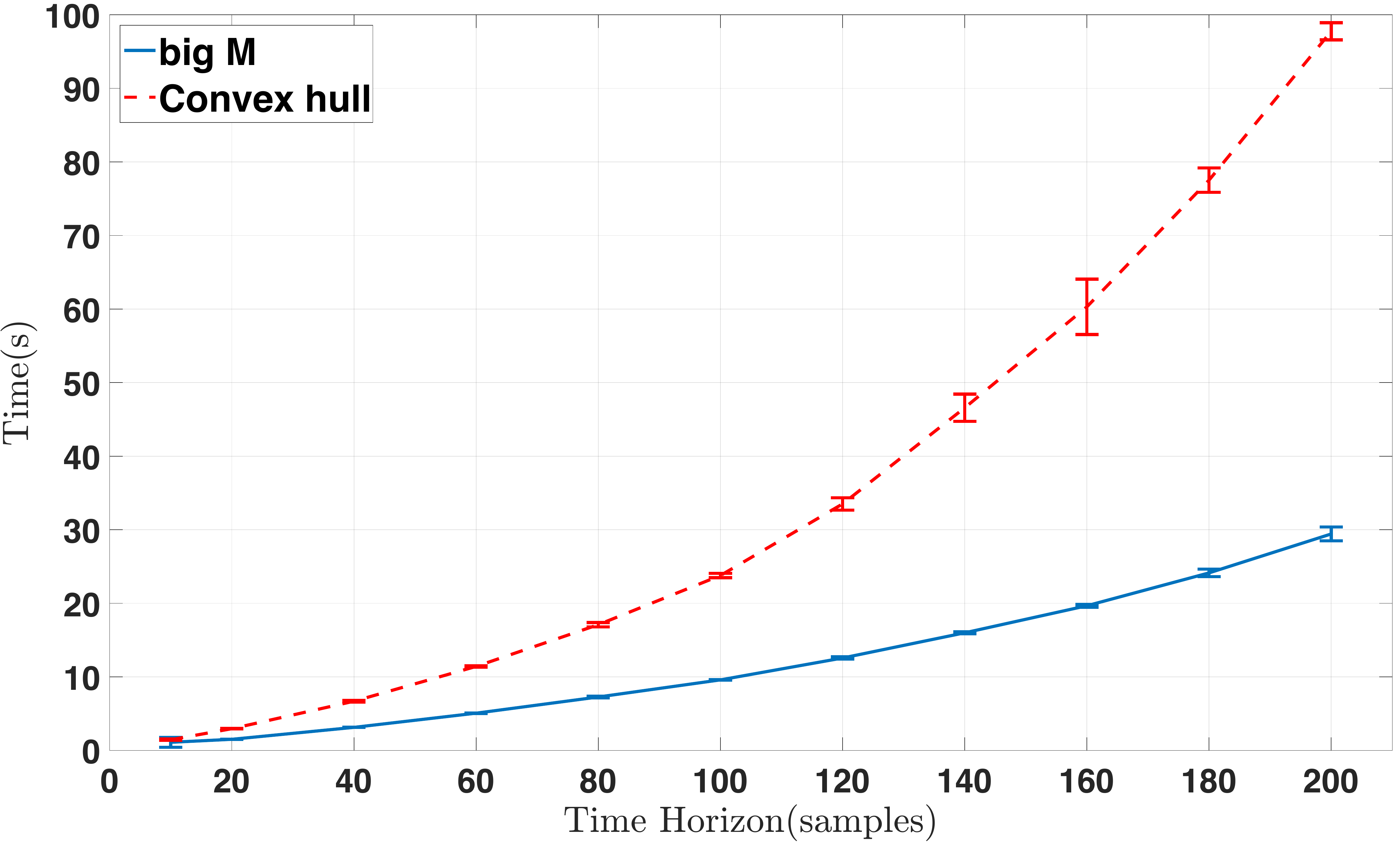}
	\caption{The average execution time for invalidation of data generated by $\G^f$ for $\G_3$ on various 
		time horizons for two formulations} 
	\label{fig:runtimecompare}
\end{figure}

Finally, to compare the complexity of convex hull formulation of model invalidation problem proposed in \cite{Harirchi2015Model}, and the big $M$ formulation that is proposed in this paper, we illustrate the average run-time for the system model with three modes, when the data is generated by faulty model for different time horizon lengths. Fig. \ref{fig:runtimecompare} shows the superiority of big $M$ formulation for the model invalidation problem.

\vspace{-0.2cm}
\subsection{Building Radiant Systems}
\subsubsection{System Model} 
\begin{figure}[h!] 
	\centering 
	\includegraphics[width=5in]{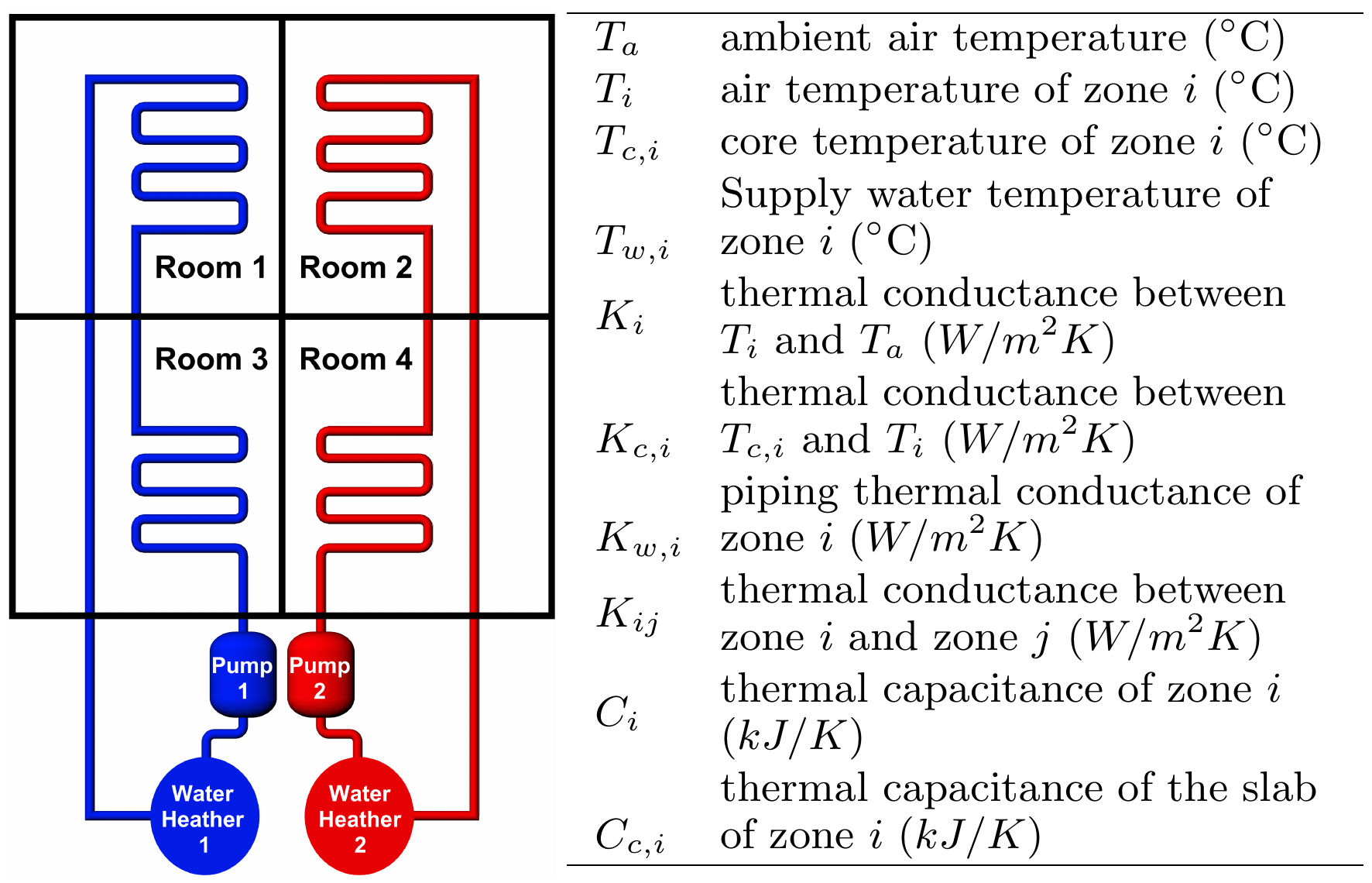} 
	\caption{Left: four zone radiant system with two pumps. Right: parameter definitions and units for the radiant system.}\label{fig:4zones+table3} 
\end{figure}

We consider a building with four rooms of the same size. The building is equipped with a radiant system with two 
pumps, adapted from \cite{NghiemEvent2013} and illustrated in Fig. \ref{fig:4zones+table3}.

There is a concrete slab (core) between the rooms and the water pipes. We consider two cores, one in contact with the water coming from pump 1 and the other in contact with water coming from pump 2.
We assume that the two pumps can either be on with known constant flow or off. 
Each pump is connected to a valve, which adjusts the constant flow of the pump. The overall system can be described by an SWA model with six states that are corresponding to the temperatures of the four rooms and the two cores. 
We assume that all of the states are measured with some measurement noise. 
The system has four modes: mode 1 corresponds to the case where both pumps are off, mode 2 is active when pump 1 is on and pump 2 is off, mode 3 represents the case when pump 1 is off and pump 2 is on, and mode 4 is when both pumps are on.

\begin{figure*}[t!] 
	\centering 
	\includegraphics[width=6in]{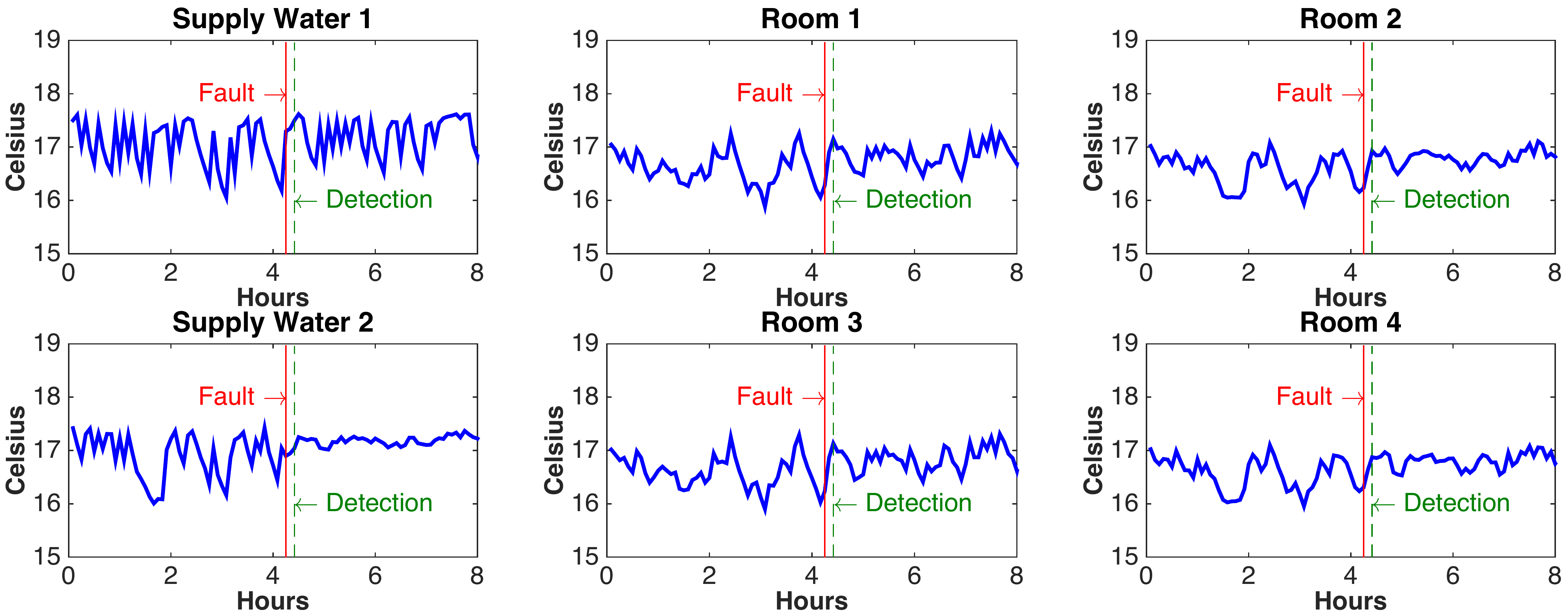} 
	\caption{Fault detection results on the data, which consists of outputs of the system until time sample $50$ and 
		the outputs of the faulty system afterwards.}\label{fig:fault} 
\end{figure*}
When both pumps are on the core temperatures, $T_{c,1}, T_{c,2}$ evolve as a function of the parameters of the system and the temperatures of the zones, $T_i, i \in \mathbb{Z}_4^+$, according to:\\
\vspace*{-1cm}
\begin{align}
C_{c,1}\dot{T}_{c,1}(t)  = & K_{c,1} (T_1-T_{c,1}) + K_{c,3}(T_3-T_{c,1}) + K_{w,1}(T_{w,1}-T_{c,1}) \label{eqn:firstpump}  \\
C_{c,2}\dot{T}_{c,2}(t)  = & K_{c,2} (T_2-T_{c,2}) + K_{c,4}(T_4-T_{c,2}) + K_{w,2}(T_{w,2}-T_{c,2}). \label{eqn:secondpump} 
\end{align}
If the first (second) pump is off, the last term in \eqref{eqn:firstpump} (in \eqref{eqn:secondpump}) becomes zero. 
The rest of state equations do not change with 
the status of the two pumps and are as follows:
\begin{align*} 
	C_{1}\dot{T}_{1}(t)  = & K_{c,1} (T_{c,1}-T_1) + K_1(T_a-T_1) + K_{12}(T_2-T_1) + K_{13}(T_3-T_1)  \\
	C_{2}\dot{T}_{2}(t)  = & K_{c,2} (T_{c,2}-T_2) + K_2(T_a-T_2) + K_{12}(T_1-T_2) + K_{24}(T_2-T_4)  \\
	C_{3}\dot{T}_{3}(t)  = & K_{c,1} (T_{c,1}-T_3) + K_3(T_a-T_3) + K_{13}(T_1-T_3) + K_{34}(T_4-T_3)  \\
	C_{4}\dot{T}_{4}(t)  = & K_{c,2} (T_{c,2}-T_4) + K_4(T_a-T_4) + K_{24}(T_2-T_4) +  K_{34}(T_3-T_4).  \\ 
	\end{align*} 
The list of parameters are given in Fig. \ref{fig:4zones+table3} and their values\footnote{$T_a = 10$, $K_1 = K_2 = \frac{1}{2.1}$, $K_3 = K_4 = \frac{1}{2.2}$, $K_{c,1} = K_{c,3} =\frac{1}{0.125}$, $K_{c,2} = K_{c,4}= \frac{1}{0.130}$, $ K_{12} = K_{13} = K_{24} = K_{34} = \frac{1}{0.16} $, $K_{w,1}= \frac{1}{0.07}$, $K_{w,2} = \frac{1}{0.05}$, 
	$ C_1 = 1900$, $C_2 = 2100$,  $C_3 = 2000$,  $C_4 = 1800$, $C_{c,1} = 3000$,  $C_{c,2} = 4000$, 
	$ T_{w,1} = T_{w,2} = 18$.} are adapted from \cite{NghiemEvent2013}.  

The discrete-time switched affine model $\G_R= (\X, \E, \U ,$ $\{ \G^\Delta_i\}_{i=1}^4)$ 
for the system is obtained with sampling time of 5 minutes, 
where $\X = \{x \mid 15 \leq x_i \leq 19 \}$, $\U= \emptyset$, $\E = \{\ETA \mid \| \ETA \| \leq 0.05 \}$, and 
$G_i$ is the discrete uncertain affine model of the $i$th mode.
We assume the uncertainty in the system is due to small changes in ambient temperature ($T_a = 10 + \delta$ with $\|\delta \| \leq 0.5$). The effect 
of the uncertainty appears in the $\f_i$ vector for each mode $i$. 
\subsubsection{Fault model} \label{subsec:fault}
We assume that the valve of the second pump stuck in middle and do not respond to on or off commands. 
In faulty mode the system has only two modes, which correspond to the ``on'' or ``off'' states of the first pump. 
This fault is modeled with a change in the second heat conductance 
parameter, which is assumed to be $K_{w,2} = \frac{0.5}{0.05}$ when the fault occurs. This indicates the fact that the heat 
transfer is cut in half, because of the slower water flow. We assume the faulty model is represented by 
$\G^f_R= (\X, \E, \U , \{ \G^{\Delta,f}_i\}_{i=1}^2)$. 
This is an incipient fault and hard to detect because the outputs can remain within the 
reasonable range, and do not change dramatically.

\subsubsection{Results}
\vspace{-0.1cm}
\paragraph{$T$-detectability:} First, we analyze the $T$-detectability of the pair 
$(\G_R, \G^f_R)$. In particular, our approach finds that for $T=8$,
this fault is detectable for the radiant system. 
The iterations from $T=1$ until detectability is verified at $T=8$ took
$36.57$ seconds with the proposed approach. 
As a comparison, satisfiability modulo theory (SMT) based approach
proposed in \cite{Harirchi2015Model} took $10008.19$ seconds.
Although finding a $T$ such that a particular fault is $T$-detectable for a specific system 
is an offline step, hence, one can tolerate a slightly slower algorithm,
the MILP based approach proposed in this paper shows significant improvement
in execution time.
%
\vspace{-0.6cm}
\paragraph{Fault detection:} Next, we demonstrate the proposed fault detection scheme on this system.
In this simulation, we generate data from the radiant system 
model $\G_R$ for $50$ samples (four hours and ten minutes), and the
persistent fault 
$\G^f_R$ becomes active at time sample $51$. 
The model invalidation problem is solved in receding horizon manner with
horizon size $8$. The fault is detected at time sample $53$. Simulation
traces are shown in Fig.~\ref{fig:fault}.
It is worth mentioning that $T$ is calculated based 
on the worst case scenario to provide guarantees for detection, but for a particular realization of input-output 
trajectory usually it is possible to detect the fault earlier as is the case in this example.
\vspace{-0.6cm}
\paragraph{Redundant sensors:}
In this subsection, we show how $T$-detectability analysis can be used for
selecting ``optimal" sensors for fault detection.
We consider five different scenarios, each corresponding to
different set of states being measured as shown in Tab. \ref{tab:Tsensors}.
For each of the scenarios, we consider the case with no uncertainty and noise in the system and fault models, e.g., $\E =\{ \mathbf{0} \}$ 
and $\Omega =\{ \mathbf{0} \}$, where $\mathbf{0}$ indicates a vector of zeros of appropriate dimension;
and also the noise and uncertainty models used in this section earlier.
We compute the minimum $T$-values for each scenario to achieve $T$-detectability
as listed in Tab. \ref{tab:Tsensors}.
From these results we see that measuring the core temperatures seems
crucial, especially in case of uncertainty and noise, for detecting this fault.

\begin{table}[tbh]
	\centering
	\caption{Value of $T$ for the five scenarios, under two model assumptions of with and without uncertainty and
		 noise.}\label{tab:Tsensors}
	\renewcommand{\arraystretch}{1.5}
	\begin{tabular}{C{3.5cm} C{3.5cm} C{3.5cm}} \hline
	Measured States	& $T$ (no uncertainty) & $T$ (with uncertainty)\\ \hline
	$T_{c,1}, T_{c,2}, T_{1}, T_{2}, T_{3}, T_{4}$ & 1 & 6 \\
	$T_{c,2}, T_{1}, T_{2}, T_{3}, T_{4}$ & 1 & 6 \\
	$T_{1},T_{2}, T_{3}, T_{4}$ & 2 & $>100$\\
	$T_{c,1},T_{c,2}, T_{1}, T_{2}$ & 2 & 6\\
	$T_{c,2},T_{1},T_{3}$ & 2 & 7\\ \hline
	\end{tabular}
\end{table}

\begin{figure}[h!] 
	\centering 
	\includegraphics[width=5in]{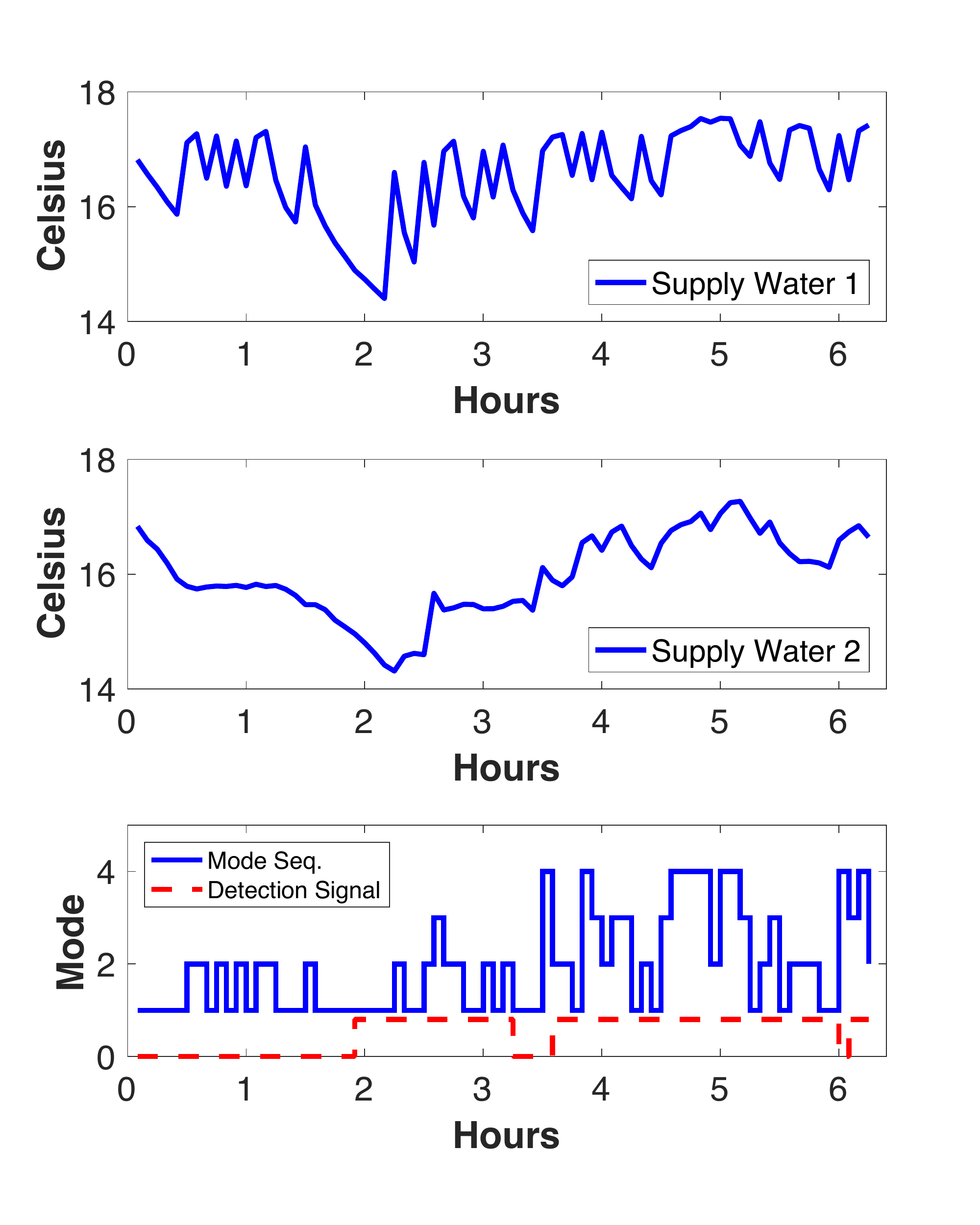} 
	\caption{Weak-detectability results, where one of the last two modes of faulty model is activated
		at sample $31$ for the first time and after sample $41$. The red dashed line in the bottom figures is $1$ when
		a fault is detected.} \label{fig:weakfault} 
\end{figure}

\paragraph{Weak-detectability:}
Consider a variant of the fault model described in Subsection \ref{subsec:fault} where the second valve can be closed, but when 
the pump is on it can only open up to the half of its capacity. Such a fault is not $T$-detectable for any finite $T$, 
because it shares two modes with the a priori system model. Problem \eqref{eqn:Tmilp} is always feasible, because 
there always exists a switching sequence that matches the fault and system model, which is the one that keeps the 
second pump always off. In order to detect such a fault, activation of the ``on'' mode for second pump is necessary, 
which corresponds to the third and fourth modes of the fault model. We capture this with an indicator of the form
$\mathcal{I}_t = (\{3,4\},1,1,=)$.
This indicator is incorporated into problem \eqref{eqn:Tmilp} via the constraint:
$\textstyle\sum_{i\in \mathbb{Z}_4^+} \sum_{j = 3}^4  d_{i,j,t-T+1} = 1$.
Solving problem \eqref{eqn:Tmilp} with this extra constraint renders $T = 8$. In order to illustrate 
the use of weak-detectability in the fault detection scheme, we generate an output sequence 
from the faulty model as follows: (i) samples 1-30 are generated by the first two modes of faulty system;
(ii) sample 31 is generated by one of the last two modes of the faulty system, and is followed by 9 samples generated 
	by the first two modes; and (iii) samples 41-75 are generated by allowing any mode of the faulty system to be active.
As illustrated in Fig. \ref{fig:weakfault}, the first detection occurs immediately after the first switch to the modes 
where second pump is ``on". The fault is also detected for all samples after $43$, where there is switching to the last 
two modes. Recall that, the mode signal is not measured by our fault detection scheme. As one can 
see the fault is not detected in the first $30$ samples, because the second pump is always ``off" during that period. 
Note that the fault does not affect the behavior of the system in the first two hours and forty minutes and detecting it is impossible.

\section{Conclusions and Discussion}\label{sec:conclusion}

In this paper, we present a fault detection scheme that 
guarantees the detection of particular faults and 
can be implemented in real-time for many applications.
The proposed scheme is applicable to an expressive class of system
and fault models, namely hidden-mode switched affine models
with parametric uncertainty. The modeling framework is further 
extended to include language constraints on the modes and
necessary and sufficient conditions for detectability using
a receding horizon scheme are provided based on
MILP. The key step in efficiency of implementation is the ability to preserve
detection guarantees while using the receding horizon scheme
as a result of $T$-detectability property. 
$T$-detectability essentially tells us how much data
is enough for detecting a fault by utilizing the knowledge of a (possibly uncertain) 
fault model and the structure in the dynamics.
With the advances in MILP-solvers, there has been
a resurge of interest in using it in control design and real-time
implementations \cite{raman2015reactive}. Here we leverage such advances in
the context of fault detection. In particular, since we only solve
MILP feasibility problems (as opposed to optimization), 
this leads to a fairly efficient
scheme as demonstrated via examples.

In the future, we are interested in applying similar ideas
to fault isolation and active fault detection. We also have
some ongoing work on fault detection for nonlinear hybrid 
systems.

\section*{Acknowledgment}
The authors would like to thank anonymous reviewers for their constructive comments and suggestions and Mustafa Kara for helpful discussions on the converse results for 
$T$-detectability. This work is supported in part by DARPA grant N66001-14-1-4045 and an 
Early Career Faculty grant from NASA's Space Technology Research Grants Program.

\bibliographystyle{unsrt}
\bibliography{IEEEabrv,modelinvalidation_abbrev,fault_abbrev}

\begin{thebibliography}{10}

\bibitem{Rajah2014Taming}
V.~Rajah.
\newblock Taming the data deluge, 2014.

\bibitem{Sznaier2014Surviving}
M.~Sznaier, O.~Camps, N.~Ozay, and C.~Lagoa.
\newblock Surviving the upcoming data deluge: A systems and control
  perspective.
\newblock In {\em IEEE CDC}, Dec 2014.

\bibitem{Weimer2013Parameter}
J.~Weimer, J.~Araujo, M.~Amoozadeh, S.~Ahmadi, H.~Sandberg, and K.~Johansson.
\newblock Parameter-invariant actuator fault diagnostics in cyber-physical
  systems with application to building automation.
\newblock In {\em Control of CPS}, pages 179--196. Springer International
  Publishing, 2013.

\bibitem{Burkart2011Nonlinear}
R.~Burkart, K.~Margellos, and J.~Lygeros.
\newblock Nonlinear control of wind turbines: An approach based on switched
  linear systems and feedback linearization.
\newblock In {\em IEEE CDC-ECC}, pages 5485--5490, 2011.

\bibitem{Sun2006Switched}
Z.~Sun.
\newblock {\em Switched linear systems: control and design}.
\newblock Springer Sc. \& Bus. Med., 2006.

\bibitem{Rosa2010Fault}
P.~Rosa, C.~Silvestre, J.~Shamma, and M.~Athans.
\newblock Fault detection and isolation of {LTV} systems using set-valued
  observers.
\newblock In {\em IEEE CDC}, pages 768--773, 2010.

\bibitem{Ozay2010Model}
N.~Ozay, M.~Sznaier, and C.~Lagoa.
\newblock Model (in)validation of switched {ARX} systems with unknown switches
  and its application to activity monitoring.
\newblock In {\em IEEE CDC}, pages 7624--7630, 2010.

\bibitem{OzayConvex2014}
N.~Ozay, M.~Sznaier, and C.~Lagoa.
\newblock Convex certificates for model (in)validation of switched affine
  systems with unknown switches.
\newblock {\em {IEEE} Trans. Autom. Control}, 59(11):2921--2932, Nov 2014.

\bibitem{cplex}
IBM~ILOG CPLEX.
\newblock User's manual for {CPLEX}.
\newblock {\em Int. Bus. Mach. Corp.}, 46(53):157, 2009.

\bibitem{Harirchi2015Model}
F.~Harirchi and N.~Ozay.
\newblock Model invalidation for switched affine systems with applications to
  fault and anomaly detection.
\newblock {\em IFAC ADHS}, 48(27):260--266, 2015.

\bibitem{Beard1971Failure}
R.~Beard.
\newblock {\em Failure accommodation in linear systems through
  self-reorganization.}
\newblock PhD thesis, MIT, 1971.

\bibitem{Jones1973Failure}
H.~Jones.
\newblock {\em Failure detection in linear systems.}
\newblock PhD thesis, MIT, 1973.

\bibitem{Simani2003Model}
S.~Simani, C.~Fantuzzi, and R.~Patton.
\newblock {\em Model-based fault diagnosis in dynamic systems using
  identification techniques}.
\newblock Springer Sc. \& Bus. Med., 2003.

\bibitem{Isermann2006Fault}
R.~Isermann.
\newblock {\em Fault-diagnosis systems: an introduction from fault detection to
  fault tolerance}.
\newblock Springer Sc. \& Bus. Med., 2006.

\bibitem{Ding2008Model}
S.~Ding.
\newblock {\em Model-based fault diagnosis techniques: design schemes,
  algorithms, and tools}.
\newblock Springer Sc. \& Bus. Med., 2008.

\bibitem{Patton2013Issues}
R.~Patton, P.~Frank, and R.~Clark.
\newblock {\em Issues of fault diagnosis for dynamic systems}.
\newblock Springer Sc. \& Bus. Med., 2013.

\bibitem{Isermann1993Fault}
R.~Isermann.
\newblock Fault diagnosis of machines via parameter estimation and knowledge
  processing--tutorial paper.
\newblock {\em Automatica}, 29(4):815--835, 1993.

\bibitem{Frank1993Advances}
P.~Frank.
\newblock Advances in observer-based fault diagnosis.
\newblock In {\em Int. Conf. on Fault Diag.: TOOLDIAG}, 1993.

\bibitem{Patton1997Observer}
R.~Patton and J.~Chen.
\newblock Observer-based fault detection and isolation: robustness and
  applications.
\newblock {\em Cont. Eng. Prac.}, 5(5):671--682, 1997.

\bibitem{Shames2011Distributed}
I.~Shames, A.~Teixeira, H.~Sandberg, and K.~Johansson.
\newblock Distributed fault detection for interconnected second-order systems.
\newblock {\em Automatica}, 47(12):2757--2764, 2011.

\bibitem{Gertler1997Fault}
J.~Gertler.
\newblock Fault detection and isolation using parity relations.
\newblock {\em Cont. Eng. Prac.}, 5(5):653--661, 1997.

\bibitem{Rosa2013Fault}
P.~Rosa and C.~Silvestre.
\newblock Fault detection and isolation of {LPV} systems using set-valued
  observers: An application to a fixed-wing aircraft.
\newblock {\em Cont. Eng. Prac.}, 21(3):242--252, 2013.

\bibitem{Nikoukhah1998Guaranteed}
R.~Nikoukhah.
\newblock Guaranteed active failure detection and isolation for linear
  dynamical systems.
\newblock {\em Automatica}, 34(11):1345--1358, 1998.

\bibitem{Nikoukhah2006Auxiliary}
R.~Nikoukhah and S.~Campbell.
\newblock Auxiliary signal design for active failure detection in uncertain
  linear systems with a priori information.
\newblock {\em Automatica}, 42(2):219--228, 2006.

\bibitem{Scott2014Input}
J.~K. Scott, R.~Findeisen, R.~D Braatz, and D.~M. Raimondo.
\newblock Input design for guaranteed fault diagnosis using zonotopes.
\newblock {\em Automatica}, 50(6):1580--1589, 2014.

\bibitem{Raimondo2016Closed}
D.~M. Raimondo, G.~R. Marseglia, R.~D. Braatz, and J.~K. Scott.
\newblock Closed-loop input design for guaranteed fault diagnosis using
  set-valued observers.
\newblock {\em Automatica}, 74:107 -- 117, 2016.

\bibitem{Garcia1997Deterministic}
E.~Garcia and P.~Frank.
\newblock Deterministic nonlinear observer-based approaches to fault diagnosis:
  a survey.
\newblock {\em Cont. Eng. Prac.}, 5(5):663--670, 1997.

\bibitem{Hammouri1999Observer}
H.~Hammouri, M.~Kinnaert, and E.~El~Yaagoubi.
\newblock Observer-based approach to fault detection and isolation for
  nonlinear systems.
\newblock {\em {IEEE} Trans. Autom. Control}, 44(10):1879--1884, 1999.

\bibitem{Pan2015Online}
W.~Pan, Y.~Yuan, H.~Sandberg, J.~Gon{\c{c}}alves, and G.~Stan.
\newblock Online fault diagnosis for nonlinear power systems.
\newblock {\em Automatica}, 55:27--36, 2015.

\bibitem{De2001Geometric}
C.~De~Persis and A.~Isidori.
\newblock A geometric approach to nonlinear fault detection and isolation.
\newblock {\em {IEEE} Trans. Autom. Control}, 46(6):853--865, 2001.

\bibitem{Mcilraith2000Hybrid}
S.~McIlraith, G.~Biswas, D.~Clancy, and V.~Gupta.
\newblock Hybrid systems diagnosis.
\newblock In {\em Hybrid Systems: Computation and Control}, pages 282--295.
  Springer, 2000.

\bibitem{Narasimhan2007Model}
S.~Narasimhan and G.~Biswas.
\newblock Model-based diagnosis of hybrid systems.
\newblock {\em {IEEE} Trans. Syst., Man, Cybern. {A}}, 37(3):348--361, 2007.

\bibitem{Deng2015Verification}
Y.~Deng, A.~D'Innocenzo, M.~Di~Benedetto, S.~Di~Gennaro, and A.~Julius.
\newblock Verification of hybrid automata diagnosability with measurement
  uncertainty.
\newblock {\em {IEEE} Trans. Autom. Control}, 61(4):982--993, 2016.

\bibitem{Grewal1976Identifiability}
M.~Grewal and K.~Glover.
\newblock Identifiability of linear and nonlinear dynamical systems.
\newblock {\em {IEEE} Trans. Autom. Control}, 21(6):833--837, 1976.

\bibitem{De2016Observability}
E.~De~Santis and M.~D. Di~Benedetto.
\newblock Observability of hybrid dynamical systems.
\newblock {\em Foundations and Trends{\textregistered} in Systems and Control},
  3(4):363--540, 2016.

\bibitem{vidal2002observability}
R.~Vidal, A.~Chiuso, and S.~Soatto.
\newblock Observability and identifiability of jump linear systems.
\newblock In {\em IEEE CDC}, volume~4, pages 3614--3619, 2002.

\bibitem{Babaali2004Observability}
M.~Babaali and M.~Egerstedt.
\newblock Observability of switched linear systems.
\newblock In {\em Int. Workshop on Hybrid Systems: Computation and Control},
  pages 48--63. Springer, 2004.

\bibitem{De2011Location}
E.~De~Santis.
\newblock On location observability notions for switching systems.
\newblock {\em Systems \& Control Letters}, 60(10):807--814, 2011.

\bibitem{Lou2009Distinguishability}
H.~Lou and P.~Si.
\newblock The distinguishability of linear control systems.
\newblock {\em Nonlinear Analysis: Hybrid Systems}, 3(1):21--38, 2009.

\bibitem{Rosa2011Distinguishability}
P.~Rosa and C.~Silvestre.
\newblock On the distinguishability of discrete linear time-invariant dynamic
  systems.
\newblock In {\em IEEE CDC-ECC}, pages 3356--3361, 2011.

\bibitem{Adnan2011Expected}
N.~Adnan, I.~Izadi, and T.~Chen.
\newblock On expected detection delays for alarm systems with deadbands and
  delay-timers.
\newblock {\em Journal of Process Control}, 21(9):1318--1331, 2011.

\bibitem{Mariton1989Detection}
M.~Mariton.
\newblock Detection delays, false alarm rates and the reconfiguration of
  control systems.
\newblock {\em International Journal of Control}, 49(3):981--992, 1989.

\bibitem{Stoorvogel2001Delays}
A.~Stoorvogel, H.~Niemann, and A.~Saberi.
\newblock Delays in fault detection and isolation.
\newblock In {\em ACC}, volume~1, pages 459--463, 2001.

\bibitem{smith1992model}
R.~Smith and J.~Doyle.
\newblock Model validation: A connection between robust control and
  identification.
\newblock {\em {IEEE} Trans. Autom. Control}, 37(7):942--952, 1992.

\bibitem{Cheng2012Convex}
Y.~Cheng, Y.~Wang, M.~Sznaier, N.~Ozay, and C.~Lagoa.
\newblock A convex optimization approach to model (in)validation of switched
  arx systems with unknown switches.
\newblock In {\em IEEE CDC}, pages 6284--6290, Dec 2012.

\bibitem{Harirchi2016Model}
F.~Harirchi, Luo Z., and N.~Ozay.
\newblock Model (in)validation and fault detection for systems with polynomial
  state-space models.
\newblock In {\em ACC}, pages 1017--1023, 2016.

\bibitem{Raman1994Modelling}
R.~Raman and I.~Grossmann.
\newblock Modelling and computational techniques for logic based integer
  programming.
\newblock {\em Comput. \& Chem. Eng.}, 18(7):563--578, 1994.

\bibitem{Bemporad1999Control}
A.~Bemporad and M.~Morari.
\newblock Control of systems integrating logic, dynamics, and constraints.
\newblock {\em Automatica}, 35(3):407--427, 1999.

\bibitem{Bertsimas2006Tractable}
D.~Bertsimas and M.~Sim.
\newblock Tractable approximations to robust conic optimization problems.
\newblock {\em Math. Prog.}, 107(1-2):5--36, 2006.

\bibitem{Egerstedt}
M.~Egerstedt and M.~Babaali.
\newblock On observability and reachability in a class of discrete-time
  switched linear systems.
\newblock In {\em ACC}, pages 1179--1180, 2005.

\bibitem{Babaali2003Pathwise}
M.~Babaali and M.~Egerstedt.
\newblock Pathwise observability and controllability are decidable.
\newblock In {\em IEEE CDC}, volume~6, pages 5771--5776. IEEE, 2003.

\bibitem{Sampath1995Diagnosability}
M.~Sampath, R.~Sengupta, S.~Lafortune, K.~Sinnamohideen, and D.~Teneketzis.
\newblock Diagnosability of discrete-event systems.
\newblock {\em {IEEE} Trans. Autom. Control}, 40(9):1555--1575, 1995.

\bibitem{Baier2008Principles}
C.~Baier and J.~Katoen.
\newblock {\em Principles of model checking}, volume 26202649.
\newblock MIT press Cambridge, 2008.

\bibitem{pasqualetti2013attack}
F.~Pasqualetti, F.~D{\"o}rfler, and F.~Bullo.
\newblock Attack detection and identification in cyber-physical systems.
\newblock {\em {IEEE} Trans. Autom. Control}, 58(11):2715--2729, 2013.

\bibitem{Shoukry2014Secure}
Y.~Shoukry, P.~Nuzzo, A.~Puggelli, A.~Sangiovanni-Vincentelli, S.~Seshia, and
  P.~Tabuada.
\newblock Secure state estimation for cyber physical systems under sensor
  attacks: a satisfiability modulo theory approach.
\newblock {\em IEEE Trans. on Autom. Control}, 2017.

\bibitem{Chong2015Observability}
M.~Chong, M.~Wakaiki, and J.~Hespanha.
\newblock Observability of linear systems under adversarial attacks.
\newblock In {\em ACC}, pages 2439--2444, 2015.

\bibitem{YALMIP}
J.~{L\"{o}fberg}.
\newblock Yalmip : A toolbox for modeling and optimization in {MATLAB}.
\newblock In {\em CACSD Conference}, Taipei, Taiwan, 2004.

\bibitem{NghiemEvent2013}
T.~Nghiem, G.~Pappas, and R.~Mangharam.
\newblock Event-based green scheduling of radiant systems in buildings.
\newblock In {\em ACC}, pages 455--460, 2013.

\bibitem{raman2015reactive}
V.~Raman, A.~Donz{\'e}, D.~Sadigh, R.~M. Murray, and S.~A. Seshia.
\newblock Reactive synthesis from signal temporal logic specifications.
\newblock In {\em Hybrid Systems: Computation and Control}, pages 239--248.
  ACM, 2015.

\end{thebibliography}

\end{document}